\documentclass[a4paper,10pt]{article}

\usepackage[a4paper,left=10mm,right=10mm,top=10mm,bottom=10mm]{geometry}
\usepackage{multirow}
\usepackage{graphicx}
\usepackage{amsfonts}
\usepackage{amsmath}
\usepackage{amssymb}
\usepackage{colortbl}
\usepackage{indentfirst}
\usepackage{booktabs}
\usepackage{verbatim}
\usepackage{color}
\usepackage{amsthm}
\usepackage{booktabs}
\usepackage{caption}
\usepackage{pdflscape}
\usepackage{algorithm}
\usepackage[noend]{algpseudocode}
\usepackage{enumitem}
\usepackage{subfigure}
\usepackage{diagbox}
\usepackage{epstopdf}

\usepackage[page,header]{appendix}

\newcommand{\BE}{{\mathbb{E}}}

\newcommand{\BN}{{\mathbb{N}}}

\newcommand{\CE}{{\mathcal{E}}}
\newcommand{\CF}{{\mathcal{F}}}

\newcommand{\CJ}{{\mathcal{J}}}

\newcommand{\CL}{{\mathcal{L}}}

\newcommand{\CN}{{\mathcal{N}}}
\newcommand{\CO}{{\mathcal{O}}}

\newcommand{\CS}{{\mathcal{S}}}

\newcommand{\beq}{\begin{equation}}
\newcommand{\eeq}{\end{equation}}
\newcommand{\bde}{\begin{definition}}
\newcommand{\ede}{\end{definition}}
\newcommand{\brm}{\begin{remark}}
\newcommand{\erm}{\end{remark}}
\newcommand{\bex}{\begin{example}}
\newcommand{\eex}{\flushright{\(\qedsymbol\)}\end{example}}
\newcommand{\bthm}{\begin{theorem}}
\newcommand{\ethm}{\end{theorem}}
\newcommand{\bprp}{\begin{proposition}}
\newcommand{\eprp}{\end{proposition}}
\newcommand{\blm}{\begin{lemma}}
\newcommand{\elm}{\end{lemma}}

\newcommand{\bprf}{\begin{proof}}
\newcommand{\eprf}{\end{proof}}
\newcommand{\ben}{\begin{eqnarray}}
\newcommand{\een}{\end{eqnarray}}

\newcommand{\dd}{\mathrm{d}}

\newcommand{\kl}{\left(}
\newcommand{\kr}{\right)}
\newcommand{\bmr}{\left(\begin{array}}
\newcommand{\emr}{\end{array}\right)}

\newcommand{\bal}{\begin{aligned}}
\newcommand{\eal}{\end{aligned}}

\newcommand{\nbr}{\nonumber}

\numberwithin{equation}{section}
\newtheorem{theorem}{Theorem}[section]
\newtheorem{lemma}[theorem]{Lemma}

\newtheorem{proposition}[theorem]{Proposition}
\newtheorem{definition}[theorem]{Definition}
\newtheorem{remark}[theorem]{Remark}

\newtheorem{claim}[theorem]{Claim}
\newtheorem{assumption}{Assumption}

\begin{document}

	\title{A synchronization-capturing multi-scale solver to the noisy integrate-and-fire neuron networks}
	
	\author{Ziyu Du \footnote{The Oden Institute for Computational Engineering and Sciences, The University of Texas at Austin, Austin, Texas 78712, USA (ziyu\_du@utexas.edu).}\  \
 Yantong Xie\footnote{School of Mathematics Science, Peking University, Beijing, 100871, China (darkoxie@pku.edu.cn).} \   \
	Zhennan Zhou\footnote{Beijing International Center for Mathematical Research, Peking University, Beijing, 100871, China (zhennan@bicmr.pku.edu.cn).}}
	\maketitle

\begin{abstract}

The noisy leaky integrate-and-fire (NLIF) model describes the voltage configurations of neuron networks with an interacting many-particles system at a microscopic level. When simulating neuron networks of large sizes, computing a coarse-grained mean-field Fokker-Planck equation solving the voltage densities of the networks at a macroscopic level practically serves as a feasible alternative in its high efficiency and credible accuracy. However, the macroscopic model fails to yield valid results of the networks when simulating considerably synchronous networks with active firing events. In this paper, we propose a multi-scale solver for the NLIF networks, which inherits the low cost of the macroscopic solver and the high reliability of the microscopic solver. For each temporal step, the multi-scale solver uses the macroscopic solver when the firing rate of the simulated network is low, while it switches to the microscopic solver when the firing rate tends to blow up. Moreover, the macroscopic and microscopic solvers are integrated with a high-precision switching algorithm to ensure the accuracy of the multi-scale solver. The validity of the multi-scale solver is analyzed from two perspectives: firstly, we provide practically sufficient conditions that guarantee the mean-field approximation of the macroscopic model and present rigorous numerical analysis on simulation errors when coupling the two solvers; secondly, the numerical performance of the multi-scale solver is validated through simulating several large neuron networks, including networks with either instantaneous or periodic input currents which prompt active firing events over a period of time.

\end{abstract}

{\small

}



\section{Introduction}

\newcommand{\ext}{\mathrm{ext}}
\newcommand{\ite}{\mathrm{int}}

To model a network of interacting neurons, the noisy leaky integrate-and-fire (NLIF) models  originated in \cite{L1907} are widely used in computational neuroscience. The NLIF models govern the depolarization potentials (abbreviated as "potentials" or "voltages" in the rest of this paper) of individual neurons in a network as functions of time and evolve them under synaptic currents within and outside the network. The action potentials of the neurons are modeled by sudden depolarizations and hyperpolarizations of the membrane potentials, meanwhile exerting synaptic currents that excite or inhibit every other neuron in the network. Specifically, we refer to an action potential as a firing event or a spike and we also say a neuron fires or spikes when it encounters an action potential. We refer to \cite{T1988,RBW2004} as reviews of the NLIF models and the recent results.

In this paper, we are concerned with an interacting neuron network consisting of \(L=L_I+L_E\) neurons, including \(L_I\) inhibitory neurons and \(L_E\) excitatory neurons. Let the membrane potential of the \(j\)-th type \(Q\) neuron at time \(t>0\) be \(V_j^Q(t)\) where \(j=1,2,\cdots,L_Q\) and \(Q\in\{E,I\}\) indicates the neuron is either excitatory or inhibitory. When a neuron doesn't receive synaptic currents from within or outside the network, its membrane potential \(V_j^Q(t)\) relaxes to a resting potential \(V_0\) exponentially fast. More specifically, the dynamics of \(V_j^Q(t)\) is modeled by the following ODE:
\beq \label{sde00}
c\frac{\dd}{\dd t}V_j^Q(t)=-g\kl V_{j}^Q-V_L\kr+I_j^{\ext,Q}(t)+I_j^{\ite,Q}(t),\,j=1,2,\cdots,L_Q,\, Q\in\{E,I\},
\eeq
where the parameter \(c\) and \(g\) are the membrane capacitance and conductance and we assume \(c=g=1\) in the rest of this paper for simplicity. The terms \(I_j^{\ite}\) and \(I_j^\ext\) model the collective currents due to interactions with other neurons in the network and from the exterior environment respectively. The external currents are modeled with a deterministic function plus a Gaussian noise:
\beq\label{ModelIext}
I_j^{\ext,Q}(t)=I_0(t)+\sigma_0\frac{\dd B_t}{\dd t},
\eeq
where \(B_t\) denotes the standard Brownian motion. When a neuron encounters a firing event, it instantaneously discharges itself. The firing events in the NLIF model are modeled as immediate shifts of the membrane potentials from a firing threshold \(V_F\)  to a resetting potential \(V_R\) i.e.
\begin{equation}\label{eq:jump}
V_j^Q(t^-)=V_F, \quad V_j^Q(t^+)=V_R,
\end{equation}
where \(V_R<V_F\). Besides, each firing event of an excitatory or inhibitory neuron exerts synaptic currents that excite or inhibit the entire network. Therefore, the external current of a neuron \(I_j^{\ite}\) is modeled as a collection of synaptic currents from every presynaptic neuron that encounters a spike: 
\beq\label{eq:int}
I_j^{\ite,Q}(t)=\sum_{s=1}^{L_I}\sum_{k\in\BN^*} J_I\delta\kl t-T^I_{s,k}\kr+ \sum_{s'=1}^{L_E}\sum_{k\in\BN^*} J_E\delta\kl t-T^E_{s',k}\kr-\sum_{k\in\BN^*}J_Q\delta\kl t-T^Q_{j,k}\kr,
\eeq
where \(\delta\) is the Dirac Delta function and the spiketime \(T_{j,k}^Q\) (\(Q\in\{I,E\}\)) record the the \(k\)-th spike of the \(j\)-th type \(Q\) presynaptic neuron in the network. The synaptic currents associated with excitatory and inhibitory presynaptic neurons are \(J_E>0\) and \(J_I<0\). We mention that quite a few studies on the integrate-and-fire neuron networks focus on many-particles NLIF model as Equation \eqref{sde00}-\eqref{eq:int} and further modifications in order to match more complicated biological phenomenon \cite{BH1999,CS2016,CR2020}. Nevertheless, in this work we consider the simulation of the entirely excitatory or inhibitory versions of the NLIF models \eqref{sde00}-\eqref{eq:int} in the simplest form, whose evolution is governed by the following interacting system:
\begin{align}
\frac{\dd}{\dd t}V_j \,\, &=-(V_{j}-V_L)+I_0(t)+\sum_{s\neq{j}}\sum_kJ\delta(t-T_{s,k})+\sigma_0\frac{\dd B_t}{\dd t},
\label{sde1} \\
\label{eq:jump1}
V_j(t^-)&=V_F, \quad V_j(t^+)=V_R,
\end{align}
where the synaptic strength \(J=J_E\) or \(J_I\).

Though simulating the NLIF model governed by the interacting system \eqref{sde00}-\eqref{eq:int} or \eqref{sde1}-\eqref{eq:jump1} at a microscopic level is believed to approximate the realistic neuron network to a great extent, the computational complexities to solve these equations directly are costly due to the large network size \(L\gg1\). To simplify the task, people seek to simulate  coarse-grained models approximating the microscopic systems and have succeed in simulating the NLIF network with more acceptable complexity \cite{AB1997,RY2013,ZZCR2013}. Among fruitful coarse-grained models, a widely-used formulation is to assume each individual neurons in the network indistinguishable and apply the mean-field approximation to describe the dynamics of an individual neuron at a macroscopic level with a mean-field voltage function \(V(t)\) \cite{AB1997,B2000,OBH2009}. The spikes are modeled by a nonhomogeneous Poisson process with intensity \(N(t)\) at time \(t\), rather than a series of discrete shifts as in Equation \eqref{sde1}. The time-dependent function \(N(t)\), named the \textbf{mean firing rate}, serves as an important characteristic in the mean-field models. Finally, the time evolution of \(V(t)\) is described a McKean-Vlosov type SDE that gives the mean-field dynamics of an individual neuron:
\beq \label{eq:MFsde0}
\begin{cases}
\dd V=\left[-(V-V_L)+I_0(t)+bN(t)\right]\dd t+\sigma_0\dd B_t,\\
V(t^-)=V_F\,\Rightarrow\, V(t^+)=V_R,\\
N(t)=\frac{\dd}{\dd t}\BE \left[\sum_{k\in\BN^*}\chi_{[0,t]}(\tau_k)\right],
\end{cases}
\eeq
where \(\tau_k\) is the \(k-\)th time that \(V(t)\) reaches the threshold \(V_F\). The synaptic currents exerted through the spikes are modeled by the term \(bN(t)\) in the drift term, where the \textbf{connectivity parameter} \(b\) describes the synaptic connectivity of the network. Naturally, \(b>0\) indicates an excitatory network, while \(b<0\) indicates an inhibitory network. Heuristically, by Ito's formula, one can derive the time evolution of the density function \(p(v,t)\) of the voltage \(V(t)\) with the following Fokker-Planck equation with a flux shift structure:
\begin{equation} \label{eq:fk0}
\begin{cases}
\partial_t p+\partial_v [ (-v+V_L+I_0(t)+bN(t)) p] - a\partial_{vv}p
=0, \quad v\in(-\infty,V_F)/\{V_R\},  \\
p(v,0)=p_0(v),\qquad p(-\infty,t)=p(V_F,t)=0, \\
p(V_R^-,t)=p(V_R^+,t),\qquad \partial_vp(V_R^-,t)=\partial_v p(V_R^+,t)+\frac{N(t)}{a},
\end{cases}
\end{equation}
where \(a=\frac{\sigma_0^2}{2}\) denotes the diffusion coefficients. The equivalence of the SDE \eqref{eq:MFsde0} and the Fokker-Planck equation \eqref{eq:fk0} is studied rigorously through an iteration scheme in the linear case, i.e. \(b=0\) \cite{LWZZ2020}. However, even a formal connection between Equation \eqref{eq:MFsde0} and \eqref{eq:fk0} is lacking when \(N(t)\) blows up, which is a possible solution behavior for either Equation \eqref{eq:MFsde0} or Equation \eqref{eq:fk0} for relatively excitatory networks.


The Fokker-Planck equation \eqref{eq:fk0} describes the microscopic voltage configuration of the network from a macroscopic level. When considering a large system of interacting neurons, simulating the Fokker-Planck equation \eqref{eq:fk0} through the efficient structure-preserving solver proposed in \cite{HLXZ2021} (namely the macroscopic solver) is of lower computational cost compared with simulating the original system \eqref{sde1}-\eqref{eq:jump1} directly (namely the microscope solver). Moreover, the macroscopic solver has already succeeded in providing credible references to the original system \eqref{sde1}-\eqref{eq:jump1} in quite a lot of cases \cite{HLXZ2021,LWXZZ2021,HHZ2022,CR2020}. However, due to the nonlinear firing mechanism, the solutions to Equation \eqref{eq:fk0} with sufficiently large connectivity parameter \(b\) may blow up at a finite time, meanwhile Equation \eqref{eq:fk0} depicts a synchronous network with profuse spikes and strong synaptic interactions. Therefore, individual neurons in these synchronous networks are not longer indistinguishable, indicating that macroscopic solver solving the mean-field Fokker-Planck Equation \eqref{eq:fk0} may not be reliable when considering synchronous networks. In order to make best use of the macroscopic models in approximating and simulating synchronous networks, people have studied some variant models of Equation \eqref{eq:fk0} incorporating additional regularization mechanisms, such as the refractory states, transmission delay or generation solutions on the time dilation approach \cite{RY2013,CS2018,DZ2022}. However, the study on how the regularized PDE models approximate the microscopic model \eqref{sde1}-\eqref{eq:jump1} is completely open, and thus the finite-time blow-up is still a formidable challenge when simulating the NLIF model with the macroscopic solver. These situations encourage us to design a novel numerical solver that not only makes best use of the macroscopic solver, but also produces reliable results when blow-up happens.

In this paper, we propose a new multi-scale scheme for the integrate-and-fire neuron networks \eqref{sde1}-\eqref{eq:jump1} that hybridizes the macroscopic solver to the Fokker-Planck equation \eqref{eq:fk0} and the microscopic solver to the original system. The core of our new scheme is to lower the computational complexity through evolving the macroscopic solver whenever the PDE model Equation \eqref{eq:fk0} provides a reliable mean-field approximation to the original system. Practically, our scheme evolves either the macroscopic solver or the microscopic solver for each temporal step, and the choice of the solvers is determined by the mean firing rate \(N(t)\): it applies the macroscopic solver when \(N(t)\) is moderately small, but switches to the microscopic solver whenever \(N(t)\) exceeds a prescribed threshold indicating that the networks become synchronous. Moreover, the macroscopic and microscopic solver are coupled together with a high-accuracy switching algorithm which transforms between the macroscopic density functions solving Equation \eqref{eq:fk0} and the microscopic voltage configurations solving the system \eqref{sde1}-\eqref{eq:jump1}.  Finally, the validity of the multi-scale solver is analyzed from two perspectives. From the perspective of numerical analysis, we show that if we expect the error of the multi-scale solver to be lower than a particular threshold, a feasible choice of the macroscopic spatial step length is \(\CO\kl L^{-\frac{1}{4}}\kr\) when simulating a network of size \(L\). This result indicates that the computational complexity of a macroscopic temporal step of our multi-scale solver is \(\CO\kl L^{\frac{1}{4}}\kr\) for an explicit discretization and \(\CO\kl L^{\frac{3}{4}}\kr\) for a semi-implicit discretization, outperforming the microscopic solver with computational complexity \(\CO(L)\) per temporal step to a great extent when \(L\gg1\). From the perspective of numerical experiments, the numerical validity and the high efficiency of our multi-scale scheme are shown through various numerical examples, including simulating synchronous networks where instantaneous or periodic input currents push the networks to an extremely synchronous status.

Finally, we remark that designing effective computational tools for large interacting neuron networks has been a central research topic in computational neuroscience.  Among the previous studies, the seminal work \cite{ZZCR2013} established a computational framework to simulate the integrate-and-fire model with a multi-scale spirit. Following this pioneering work, a few subsequent works \cite{ZSRT2019,LL2022,SZT2020} are devoted to capturing the synchronization behavior with augmented systems or data-driven approaches. However, to the best of our knowledge, the existing methods rely on either formal scientific arguments or machine-learning methods in certain key steps of the whole algorithm, which do not yield rigorous analysis. With the mean-field approximation as the theoretical foundation, our work considers a multi-scale scheme integrating the macroscopic solver and the microscopic solver with an indicator of the network activity, and the switching only takes place before the network becomes synchronous or when it has completely become homogeneous.  Such designs allow quantitative analysis of our algorithm and naturally imply practical strategy in choosing simulation parameters, and these features distinguish our work from the former ones.

The rest of this paper is outlined as follows. In Section 2, we elaborate the mean-field approximation and the involved microscopic and macroscopic model describing the leaky-and-fire neuron network, which serve as the theoretical foundation of our multi-scale solver. In Section 3, we present our multi-scale solver and the high-precision switching algorithm coupling the macroscopic and microscopic solver. In section 4, we establish rigorous numerical analysis on the simulation error and give practical strategy in choosing simulation parameters. Finally, in Section 5, we firstly present numerical experiments that guarantee the numerical reliability between the macroscopic solver and microscopic solver. Then the numerical performance of the multi-scale solver is validated through simulating synchronous networks.

\section{The mean-field approximation and the mean-field limits of the NLIF model}

The aim of this paper is to design a numerical solver to the dynamical system \eqref{sde1}-\eqref{eq:jump1} that makes the best of the corresponding mean-field Fokker-Planck equation to enhance computation efficiency. In this section, we clarify homogeneous assumptions that are sufficient to guarantee the mean-field limit and explore the cases where the mean-field approximation may break down. The analysis concerning the mean-field approximation constitutes the theoretical foundation of the multi-scale solver to be established in later sections.  Without special instruction, we only consider purely excitatory or inhibitory networks.



\subsection{The mean-field approximation and the homogeneous assumptions}

In this subsection, we introduce the mean-field stochastic model as an approximation to the interacting particle system \eqref{sde1} -\eqref{eq:jump1} and  validate the approximation with certain homogeneous assumptions. 


In  the mean-field regime, the neurons within the network are assumed to be almost indistinguishable and loosely interacting. Instead of tracing the voltages of each neuron, we focus on the potential of a representative neuron, denoted by \(V(t)\), 
which is governed by the following SDE:
\beq\label{EqMF0}
\frac{\dd V}{\dd t}=-(V-V_L)+I_\infty^\ext(t)+ I_\infty^\ite(t),
\eeq 
where the term \(I_\infty^\ext(t)\) and \(I_\infty^\ite(t)\) describe the mean-field external current and the collective synaptic current respectively. The external current \(I_\infty^\ext(t)\) is modeled as the same deterministic process plus a Gaussian noise as Equation \eqref{ModelIext}.

The modeling of \(I_\infty^\ite(t)\) follows two steps: firstly, we adopt the diffusion approximation to the collective synaptic current \(I^\ite(t)\), which is denoted by \(\bar{I}^\ite(t)\); secondly, we consider the mean-field of the diffusion approximation \(\bar{I}^\ite(t)\) by taking the network size \(L\to\infty\). The core of the modeling is to parameterize the spike times of each neuron with a time-dependent mean firing rate function \(N(t)\):

\begin{assumption}[Mean firing rate]\label{HA1}
The spike times of each neuron in the network are assumed to be statistically uncorrelated, which follow a  nonhomogeneous Poisson process. More precisely, the probability of each neuron emitting a spike in an infinitesimal time \(\dd t\) is given by \(N(t)\dd t\). The intensity function \(N(t)\) is termed the \textbf{mean firing rate}.
\end{assumption}

With the mean firing rate \(N(t)\) prescribed, the mean and variance of the synaptic currents of each neuron are given by \(LJN(t)\dd t\) and \(LJ^2N(t)\dd t\) in an infinitesimal time \(\dd t\). Adopting the diffusion approximation, we approximate the collective point process \(I^\ite(t)\) with a smooth diffusion process \(\bar{I}^\ite(t)\) that shares the same mean and variance with \(I^\ite(t)\) \cite{AB1997,B2000,BH1999}:
\beq \label{PoissonDes}
I^\ite(t)\dd t=\kl\sum_{j'\neq{j}}\sum_kJ\delta(t-T_{j',k})\kr\dd t\approx\bar{I}^\ite(t)\dd t= LJN(t)\dd t+\sqrt{LJ^2N(t)}\dd B_t.
\eeq
We refer to Appendix A for a rigorous derivation of the diffusion approximation \eqref{PoissonDes} through the central limit theorem.


We remark that the diffusion approximation \(\bar{I}^\ite(t)\) still models the synaptic interactions of a finite-neuron network rather than the mean-field representative potential \(V(t)\). 
In this limit \(L\to\infty\)., we assume that the synaptic strength depends on the network size, i.e. \(J=J(L)\) with a proper choice of scaling. The following assumption is made to ensure that synaptic currents \(\bar{I}^\ite(t)\) of a \(L-\)neurons network has a well-defined mean-field limit \(I_\infty^\ite(t)\):

\begin{assumption}[Mean-field limit]\label{AssHom}
Consider a family of neuron networks with size \(L\), where each \(L-\)neuron network is governed by the  system \eqref{sde1}-\eqref{eq:jump1}. We assume that the synaptic strengths \(J=J(L)\) and within the time interval \([T_1, T_2)\) the firing rate \(N(t)\) satisfies :
\beq\label{HomoAssu}
\sup_{t\in[T_1,T_2)}\left[\sup_{L\in\BN^*}\Bigl( L|J| N(t)+N(t)\Bigr)\right]<\infty.
\eeq
\end{assumption}

Assumption 1 and Assumption 2 are  called \textbf{the homogeneous assumptions} in the rest of this paper. We highlight that Assumption 2 not only requires that the firing rate \(N(t)\) is bounded in the time interval \([T_1, T_2)\), but also requires that  the mean of total synaptic strength \(\bar{I}^\ite(t)\) is uniformly bounded as \(L\to+\infty\), i.e.  \(\sup_{L\in\BN^*} L|J|N(t)<\infty\) in \(t\in[T_1,T_2)\). Heuristically, the assumption implies that  the \(L-\)neuron networks are weakly correlated so that they can be approximated with a mean-field network. 

The choice of scaling \(J=J(L)\) that ensures Assumption \ref{AssHom} is obviously not unique and the most widely-used choice is to let \(\lim_{L\to+\infty}LJ(L)=b\) \cite{AB1997,BH1999,B2000}, where \(b\) is a constant named \textbf{the connectivity parameter}. The connectivity parameter models the strength of interacting of the network: \(b>0\) describes excitatory networks and \(b<0\) inhibitory networks. We may assume \(\lim_{L\to+\infty}LJ(L)=b\) without special instructions in the rest of the paper, then the Assumption 2 is satisfied as long as the firing rate is bounded. We further yield the mean-field limit variance \(\lim_{L\to\infty}L\kl J(L)\kr^2=0\), then the mean-field synaptic current \(I_\infty^\ite(t)\) is modeled as follows
\beq \label{MFIint}
I_\infty^\ite(t)\dd t=bN(t)\dd t,
\eeq
where the diffusion term vanishes as \(L\to+\infty\).

\begin{remark}
Here we mention that we may establish the diffusion approximation to excitatory-inhibitory networks \eqref{sde00}-\eqref{eq:int} as shown in \cite{CS2016}, where the authors define two firing rates \(N_E(t)\) and \(N_I(t)\) indicating the probability to emit a spike of an excitatory/inhibitory neuron:
\beq 
\bar{I}^\ite(t)\dd t= \kl L_EJ_EN_E(t)-L_IJ_IN_I(t)\kr\dd t+\sqrt{L_EJ_E^2N_E(t)+L_IJ_I^2N_I(t)}\dd B_t.
\eeq
However, the mean-field limit from \(\bar{I}^\ite(t)\) to \(\bar{I}_\infty^\ite(t)\) in this coupled network case and the corresponding homogeneous assumptions are quite complicated due to the influence from two independent firing rates. More specifically, since the drift and diffusion of the collective currents \(\bar{I}^\ite(t)\) may be of the same order, the diffusion term of the mean-field synaptic currents \(\bar{I}_\infty^\ite(t)\) may not vanish as \(L\to\infty\).
\end{remark}

Combine Equation \eqref{MFIint} with Equation \eqref{EqMF0}, we arrive at the mean-field dynamics of the membrane potential \(V(t)\):
\beq\label{eq:MFsde}
\dd V=\left[-(V-V_L)+I_0(t)+bN(t)\right]\dd t+\sigma_0\dd B_t.
\eeq
When emitting a spike, the potential \(V(t)\) discharges itself and instantaneously resets to \(V_R\), i.e.
\beq\label{FiringStr}
V(t^-)=V_F\quad\Rightarrow\quad V(t^+)=V_R.
\eeq
The firing rate \(N(t)\) is determined self-consistently through \(V(t)\) through considering the counting process of the spiketimes:
\beq
M(t)=\sum_{k\in\BN^*}\chi_{[0,t]}(\tau_k),
\eeq
where \(\chi\) stands for the indicator function. Here the hitting time \(\tau_k\) denotes the sequence of hitting time, i.e.
\beq\label{tauk}
\tau_k=\inf\left\{t\geq \tau_{k-1}:V(t)=V_F\right\},\, k\in\BN^*,
\eeq
where we assume \(\tau_0=0\) without confusion. Then \(N(t)\) is given by the derivative of expectation of the counting process \(M(t)\):
\beq\label{eq:Ntdecide}
N(t)=\frac{\dd}{\dd t}\BE M(t),
\eeq
where \(\BE\) denotes the expectation. With SDE \eqref{eq:MFsde} equipped with Equation \eqref{eq:Ntdecide}, we arrive at a self-consist system solving the mean-field voltage function \(V(t)\) and the mean-firing rate \(N(t)\). In this paper, the SDE \eqref{eq:MFsde} serves as the foundation of the coarse-grained approximations to the neuron networks. Finally, we refer to \cite{DIRT2015,DIRT2015-2} as the theoretical study of the Mckean-Vlasov type SDE \eqref{eq:MFsde}.



\subsection{The Fokker-Planck equations}

A number of studies of the integrate-and-fire neuron networks are focused on the Fokker-Planck equations, either related to the many-particles system \eqref{sde1}-\eqref{eq:jump1} \cite{AB1997,RY2013} or the mean-field SDE \eqref{eq:MFsde} \cite{CCP2011,CGGS2013,RS2021}. The Fokker-Planck equations corresponding to \(L-\)neuron networks solve the joint density \(p(V_1,\cdots,V_L,t)\) of the \(L\) neurons, which yields a \(L\)-dimensional  PDE. Though we may apply methods of separation of variables to lower the number of dimensions, the reduced models often lose certain synchronous properties of the network \cite{RY2013}. In this subsection, we study the Fokker-Planck equation derived by the mean-field SDE \eqref{eq:MFsde}, which is an one-dimensional PDE.

We refer to the density function \(p(v,t)\) as the probability of finding a neuron of voltage \(v\) at a given time \(t\) in the network. Heuristically, by Ito's formula, one can derive the time evolution of \(p(v,t)\) through the mean-field SDE \eqref{eq:MFsde}. The resulting PDE is the following Fokker-Planck equation with a flux shift:
\begin{equation} \label{eq:fk}
\begin{cases}
\partial_t p+\partial_v [ (-v+V_L+I_0(t)+bN(t)) p] - a\partial_{vv}p
=0, \quad v\in(-\infty,V_F)/\{V_R\},  \\
p(v,0)=p_0(v),\qquad p(-\infty,t)=p(V_F,t)=0, \\
p(V_R^-,t)=p(V_R^+,t),\qquad \partial_vp(V_R^-,t)=\partial_v p(V_R^+,t)+\frac{N(t)}{a},
\end{cases}
\end{equation}
where \(a=\frac{\sigma_0^2}{2}\) denotes the diffusion coefficients. At the macroscopic level, the synaptic effect of the firing events is modeled by the term \(bN(t)\) in the drift in Equation \eqref{eq:fk}. The initial condition $p_0(v)$ satisfies
\beq
\int_{-\infty}^{V_F} p_0(v) \dd v=1.
\eeq
In order to make \(p(v,t)\) a probability density function for any \(t>0\), the firing rate \(N(t)\) is chosen as follows
\begin{equation} \label{eq:Nt}
N(t)=-a\partial_v p(V_F,t) \ge 0.
\end{equation}

Due to the nonlinear firing mechanism, the Fokker-Planck equation \eqref{eq:fk0} produces diverse properties according to different choices of the connectivity parameters \(b\) and the complicated solution phenomena yield difficulties to the macroscopic solver. For \(b\leq 0\) and \(b>0\) relatively small, Equation \eqref{eq:fk} describes  networks with weak firing events and its connection with the original system \eqref{sde1}-\eqref{eq:jump1} has been studied in \cite{DIRT2015,DIRT2015-2,CR2020}. For \(b>0\) sufficiently large, Equation \eqref{eq:fk} admits solutions with strong firing events, where \(N(t)\) blows up at a finite time. In other words, when finite-time blow-ups occur, the networks become considerably synchronous  with profuse spikes and frequent synaptic interactions between neurons that violate the indistinguishability assumptions of the macroscopic models \cite{CCP2011,HLXZ2021}, hence Equation \eqref{eq:fk} no longer approximates the microscopic system \eqref{sde1}-\eqref{eq:jump1} (see Section \ref{sec:mfe} for more discussions).




To conclude this subsection, we give a formal derivation relating Equation \eqref{eq:fk} and the SDE \eqref{eq:MFsde} following the spirit of the iteration scheme proposed in \cite{LWZZ2020}, where it is proved  that the probability density of the potential function \(V(t)\) solving the SDE \eqref{eq:MFsde} is a classical solution to Equation \eqref{eq:fk} for the special case \(b=0\). We decompose the solution \(p(v,t)\) to Equation \eqref{eq:fk} as a summation of sub-density functions \(p_k(v,t)\),
\beq
p(v,t)=\sum_{k=0}^{+\infty}p_k(v,t),
\eeq
where \(p_k(v,t)\) represents the probability of finding a neuron with voltage \(v\) that has already spiked \(k\) times at time \(t\). Each sub-density \(p_k(v,t)\) matches a Fokker-Planck type sub-PDE (we assume \(V_L=0\), \(a=1\) and \(I_0(t)=0\) simplicity): for \(k\geq 1\) we have
\begin{equation} \label{eq:fksub}
\begin{cases}
\partial_t p_k+\partial_v [(-v+bN(t)) p_k] - \partial_{vv}p_k
=0, \quad v\in(-\infty,V_F)/\{V_R\},  \\
p_k(V_R^-,t)=p_k(V_R^+,t),\qquad \partial_vp_k(V_R^-,t)=\partial_v p_k(V_R^+,t)+f_{\tau_k}(t),
\end{cases}
\eeq
where the source of the flux shift  is \(f_{\tau_{k}}(t)\) denoting the density function of the \(k\)-th spiking time \(\tau_{k}\) (see Equation \eqref{tauk} for definition of \(\tau_k\)). For \(k=0\), the sub-density matches a Fokker-Planck type sub-PDE without flux source at \(V_R\):
\begin{equation} \label{eq:fksub0}
\begin{cases}
\partial_t p_0+\partial_v [(-v+bN(t)) p_0] - \partial_{vv}p_0
=0, \quad v\in(-\infty,V_F)/\{V_R\},  \\
p_0(V_R^-,t)=p_0(V_R^+,t),\qquad \partial_vp_0(V_R^-,t)=\partial_v p_0(V_R^+,t).
\end{cases}
\eeq
The spike time density \(f_{\tau_{k}}\) is determined by the \(k-1\)-th sub-density \(p_{k-1}(v,t)\) via
\beq
f_{\tau_{k-1}}(t)=-\frac{\partial f_{k-1}}{\partial x}(t,V_F).
\eeq
Therefore, in the iteration scheme, the escaping boundary flux of the \(k-1-\)th sub-density \(f_{k-1}(v, t)\) serves as the singular flux source for \(f_k(v, t)\). Note that the firing rate \(N(t)\) counts for the firing events all over the network, and thus \(N(t)\) is the summation of the spike time densities:
\beq
N(t)=\sum_{k=0}^\infty  f_{\tau_{k-1}}(t).
\eeq
Clearly, the firing rate $N(t)$ couples all sub-PDE's together through the drift term. 
 


\subsection{Beyond the homogeneous assumptions} \label{sec:mfe}


In previous sections, we provide a formal derivation of the SDE \eqref{eq:MFsde} and Fokker-Planck equation \eqref{eq:fk} as mean-field approximations to the interacting system \eqref{sde1}-\eqref{eq:jump1} under the homogeneous assumptions. However, the mean-field models fail to represent the network when the network becomes overly synchronous and the homogeneous assumptions break down. In this section, we aim to study the interacting network in these cases beyond the homogeneous assumptions.

Consider an interacting network of \(L\) neurons governed by Equation \eqref{sde1}-\eqref{eq:jump1}. If the network is excitatory and fires frequently, it is possible that the synaptic current exerted by one spike causes  other neurons to spike at the same time and trigger a chain reaction. We refer to the firing events that involve more than one spiking neurons as \textbf{multiple firing events (MFEs)}. More specifically, if an MFE involves a large portion of neurons in the network, the network becomes strongly correlated transiently which likely violates the homogeneous assumptions proposed previously. From a mathematical perspective, the existence of the MFEs indicates that the interacting system \eqref{sde1}-\eqref{eq:jump1} is actually not well-posed when a group of neurons fire in a chain reaction. As a result, we need to regulate the order in which these subsequent firing events take place when an MFE is induced. Following the spirit of \cite{RY2013,DIRT2015-2,CR2020}, we model the MFEs with the so-called \textbf{cascade mechanism} for excitatory networks. Firstly, suppose the subset
\beq
\Gamma_0=\left\{j:V_j(t^-)=V_F\right\},
\eeq
is a non-empty set at time \(t\). At this time, we establish a "cascade time axis" to record the spiking neurons under the synaptic effects at the same time \(t\). For example, the presynaptic neurons in \(\Gamma_0\) release synaptic currents of size \(J|\Gamma_0|\) on the rest of the neurons and trigger the neurons whose voltages are in the interval \(\left[V_F-J|\Gamma_0|,V_F\right]\) reach the firing threshold \(V_F\). In other words, the neurons in the subset
\beq
\Gamma_1=\left\{j\notin\Gamma_0:V_j(t^-)+J|\Gamma_0|\geq V_F\right\},
\eeq
fire due to the synaptic strengths from neurons in \(\Gamma_0\). Iteratively, for every \(k\in\BN^*\), we may define the \(k\)-th cascade subset as the neurons that reach the firing threshold under the synaptic kicks from \(\cup_{j=0}^{k-1}\Gamma_j\):
\beq
\Gamma_k=\left\{j\notin\cup_{l=0}^{k-1}\Gamma_l:V_j(t^-)+J\left|\cup_{l=0}^{k-1}\Gamma_l\right|\geq V_F\right\}.
\eeq
The cascade continues until \(\Gamma_K=\emptyset\) for some \(K\geq 1\). At this time, the cascade ceases and we define
\beq
\Gamma=\cup_{l=0}^K\Gamma_l,
\eeq
as the spiking neurons in the MFE. We now fold the cascade axis. As an immediate aftermath of an MFE, the membrane potentials of the neurons are updated as the aggregated effects of the spiking neurons involved in the cascade axis as follows:
\beq \label{MFErule1}
V_j(t^+)=\begin{cases}
V_j(t^-)+J|\Gamma|&,j\notin\Gamma,\\
V_j(t^-)+J|\Gamma|-\kl V_F-V_R\kr&,j\in\Gamma.
\end{cases}
\eeq
We remark that although the MFE model with the updating rule above yields more mathematical understanding in the mean-field limit \cite{DIRT2015-2}, it is not entirely physical. In fact, a spiking neuron enters a transient refractory state, during which it does not react to any synaptic currents \cite{RY2013}. Thus, we can modify the MFE updating rule by
\beq \label{MFErule2}
V_j(t^+)=\begin{cases}
V_j(t^-)+J|\Gamma|&,j\notin\Gamma,\\
V_R &,j\in\Gamma.
\end{cases}
\eeq
Both rules \eqref{MFErule1} \eqref{MFErule2} will be considered in the algorithms. 

Finally, we summarize the concepts concerning the MFE and cascade mechanism in the following definition:
\bde  
For excitatory network \eqref{sde1}-\eqref{eq:jump1} where \(J>0\), an MFE occurs if and only if the subset \(\Gamma_0\) and \(\Gamma_1\) are both non-empty and the neurons contained in \(\Gamma\) are neurons involved in this MFE. The quantity of neurons involved \(|\Gamma|\) is called the \textbf{MFE size} and the ratio \(\alpha=\frac{|\Gamma|}{L}\) is called the \textbf{MFE proportion}.
\ede

Now that we have defined MFEs in a \(L-\)neurons network governed by \eqref{sde1}-\eqref{eq:jump1}, it is natural to ask how to depict the MFE in the mean-field models where the network size \(L\to\infty\). Though it is straightforward to think that the dynamics of neurons may not be homogeneous with MFE onset, a quantitative  condition determining whether an MFE violates the homogeneous assumptions is lacking.  Nevertheless, we give a special case where the MFEs sufficiently break the homogeneous assumptions in the following.



Recall that we have assumed the synaptic strengths depend on the network size in the mean-field limit, i.e \(J=J(L)\). Now we explore the breakdown of the mean-field approximation via studying the following alarming case. 

\begin{claim}
The homogeneous assumptions are violated if an MFE occurs at time $t$ and \(\alpha(L)\) converges to a constant \(\alpha\in(0,1]\) as \( L \to \infty \).
\end{claim}

We show the proof by contradiction. Assume the homogeneous assumptions are satisfied in the case. Different from the single firing events, during the MFE each neuron in the network receives synaptic kicks of size \(|\Gamma|J\) . Therefore the mean of the synaptic currents is given by \(L|\Gamma(L)|J(L)N(t)\dd t\) in an infinitesimal time \(\dd t\), where \(N(t)\) is the mean firing rate function. Take \(L\to\infty\), with the condition \(\lim_{L\to\infty}LJ(L)=b\), we arrive at
\beq
\lim_{L\to\infty}L|\Gamma(L)|J(L)N(t)=bN(t)\lim_{L\to\infty}L\alpha(L)=\infty,
\eeq
which contradicts Assumption 2. In other words, the fact that the MFE proportion doesn't converge to 0 indicates that the synchronous firing events occur, hence the homogeneous assumptions surely fail.

\brm
We highlight that the notation of the cascade axis indicates that each neuron can not fire more than one spike in a MFE. With  the updating rule \eqref{MFErule1}, this condition holds if we assume \(LJ<V_F-V_R\):
\beq 
V_j(t^+)=V_j(t^-)+J|\Gamma|-\kl V_F-V_R\kr\leq V_j(t^-)+LJ-\kl V_F-V_R\kr< V_j(t^-)\leq V_F.
\eeq
For a highly excitatory network with \(J\geq\frac{V_F-V_R}{L}\), the neurons have a probability to spike more than once in the cascade axis, which violates the definition of the MFE, see \cite{DIRT2015-2}. However,  a single neuron cannot perform multiple spikes during an MFE due to the refractory state, and the modified MFE updating rule \eqref{MFErule2} naturally excludes such phenomenon  \cite{RY2013}. 
\erm


\section{Numerical schemes}


To simulate the integrate-and-fire model, two feasible choices are to simulate the interacting system \eqref{sde1}-\eqref{eq:jump1} directly (referred to as \textbf{the microscopic scheme}) and to simulate the mean-field Fokker-Planck equation \eqref{eq:fk} (referred to as \textbf{the macroscopic scheme}) as an approximation to the original model. Solving the PDE \eqref{eq:fk} with the structure-preserving finite-difference scheme originated in \cite{HLXZ2021} and inherited in \cite{HHZ2022} has the advantages of low computational complexity and high  accuracy, but suffers from the risk of producing unreliable results when the network becomes synchronous and the homogeneous assumptions fail. In this section, we propose a multi-scale solver to \eqref{sde1}-\eqref{eq:jump1} combining the macroscopic scheme with the microscopic scheme. For each temporal step, the multi-scale scheme evolves either the macroscopic solver or the microscopic solver, determined by the firing activity of the system. To elaborate our scheme, we first introduce the one-step temporal iteration of the macroscopic scheme and the microscopic scheme respectively, and then consider the switching scheme between the two schemes. Finally, we present an error analysis of the switching scheme in the following section.

\subsection{The macroscopic scheme}

We begin with a brief review to the macroscopic scheme and the involved Scharfetter-Gummel symmetric reconstruction for the convection-diffusion operator. The scheme is positivity-preserving and preserves the relative entropy structure in the semi-discrete form when \(b=0\). We refer to \cite{HLXZ2021,HHZ2022} for detailed numerical analysis on the scheme. 

Firstly, we introduce the computation mesh. Assume \(p(v,t)\) vanishes fast enough as \(v\to-\infty\) and we truncate the domain into the bounded region \([V_{\min},V_F]\times [T_{\min},T_{\max}]\). We then divide \([V_{\min},V_F]\) and \([T_{\min},T_{\max}]\) into \(n_s\) and \(n_t\) sub-intervals with size \(\Delta v=\frac{V_F-V_{\min}}{n_v}\) and \(\Delta t=\frac{T_{\max}-T_{\min}}{n_t}\). So the grid points are assigned as follows
\beq 
\begin{cases}v_i=V_{\min }+i \Delta v & i=0,1,2, \cdots, n_v, \\  t_k=T_{\min}+k \Delta t & k=0,1,2, \cdots, n_t.\end{cases}
\eeq
Then \(p_{i,k}\) and \(N_k\) represents the numerical approximation of \(p(v_i,t_k)\) (\(i\neq 0,n_v\)) and \(N(t_k)\). The numerical density \(p_{0,k}\) and \(p_{n_v,k}\) at the boundary are forced to be 0. Since \(N(t)=-a\left.\frac{\partial}{\partial v}p(v,t)\right|_{v=V_F}\), we approximate the derivative with the first-order finite difference
\beq\label{NkDeter}
N_k=-\frac{ap_{n_v-1,k}}{\Delta v},
\eeq
where we have applied the boundary condition \(p_{n_v,k}=0\). Besides, we always assume that there is an index \(r\) satisfying \(v_r=V_R\), which  treats the derivative’s discontinuity at \(v=V_R\). 

Secondly, we apply the following Scharfetter-Gummel reformulation on Equation \eqref{eq:fk}, which combines the drift term and the diffusion term:
\begin{equation} \label{eq:SG1}
\partial_t p(v,t)-a\partial_v\left(M(v,t)\partial_v\left(\frac{p(v,t)}{M(v,t)}\right)\right)=0,
\end{equation}
where
\beq
\label{eq:M}
M(v,t)=\exp\left(-\frac{(v-V_L-I_0(t)-bN(t))^2}{2a}\right).
\end{equation}
Then the macroscopic scheme adopts the spatial center difference discretization to Equation \eqref{eq:SG1}:
\begin{equation}
\label{sch:modflux}
\frac{p_{i,k+1}-p_{i,k}}{\Delta t}+ \frac{F_{i+\frac{1}{2},k+1}-F_{i-\frac{1}{2},k+1}}{\Delta v}=0,
\end{equation}
where the numerical flux at the semi-grid points \(v_{i+\frac{1}{2}}=\frac{v_i+v_{i+1}}{2}\) is approximated semi-implicitly:
\beq\label{sch:fluxdis}
F_{i+\frac{1}{2},k}=\begin{cases} -a\frac{M_{i+\frac{1}{2}}^{N_{k}}}{\Delta v}\left(\frac{p_{i+1,k+1}}{M_{i+1}^{N_{k}}}-\frac{p_{i,k+1}}{M_{i}^{N_{k}}}\right)-N_{k+1}\eta\left(v_{i+\frac{1}{2}}-V_R\right)&,\,i=0,1,\cdots,n_v-2\\
0&,\,i=-1,n_v-1,
\end{cases}
\eeq
where
\beq 
M_{i}^{N_{k}}=\exp\left(-\frac{\left[v_i-V_L-I_0(t_k)-bN_{k}\right]^2}{2a}\right),\,M_{i+\frac{1}{2}}^{N_{k}}=\frac{2M_{i+1}^{N_{k}}M_{i}^{N_{k}}}{M_{i+1}^{N_{k}}+M_{i}^{N_{k}}}.
\eeq
Here \(\eta(\cdot)\) denotes the Heaviside function and the term \(N_k\eta\left(v_{i+\frac{1}{2}}-V_R\right)\) counts for the flux shift from \(V_F\) to \(V_R\). Note that by the finite volume construction, $p_{i,k}$ is interpreted at the $i-$th cell average of the density function at time $t_k$.  

Thus, we have completely constructed a semi-implicit scheme for the Fokker-Planck equation. If we directly apply an explicit treatment to the fluxes \eqref{sch:fluxdis}, we instead obtain a fully explicitly scheme. Without special instructions, we use the semi-implicit scheme when adopting the macroscopic solver in the rest of this paper. We denote the one-step temporal iteration of the macroscopic scheme by
\beq \label{PDEscheme}
\kl \{p_{i,k+1}\}_{i=1}^{n_v},N_{k+1}\kr=\CS_{\text{mac}}\kl \{p_{i,k}\}_{i=1}^{n_v},N_{k}\kr,
\eeq
where \(\CS_{\text{PDE}}\) denotes the discrete time evolution rule.

\subsection{The microscopic scheme}

In contrast to the macroscopic scheme, the microscopic scheme directly simulates the system \eqref{sde1}-\eqref{eq:jump1}. We evolve the system  with the Euler-Maruyama scheme with an implementation of the cascade mechanism for possible MFEs.

Firstly, we partition the temporal computing domain \([T_{\min},T_{\max}]\) into \(n_t\) sub-intervals with size \(\Delta t=\frac{T_{\max}-T_{\min}}{n_t}\). So the grid points follows \(t_k=T_{\min}+k\Delta t\) with \(k=0,1,\cdots,n_t\). \(V_{j,k}\) and \(M_k\) denotes the numerical approximation of the voltage \(V_j(t_k)\) and the times of spikes that occur during time \([T_{\min},t_k]\) respectively. Naturally, we have \(M_0=0\). For each temporal step, we update the voltage configuration of the \(L\) neurons \(\left\{V_{j,k}\right\}_{j=1}^L\) and the counting process \(M_k\) recording the number of firing events occur in the network:
\beq 
\kl \left\{V_{j,k+1}\right\}_{j=1}^L,M_{k+1}\kr=\CS_{\text{mic}}\kl \left\{V_{j,k}\right\}_{j=1}^L,M_k\kr.
\eeq
 
The evolution of \(\left\{V_{j,k}\right\}_{j=1}^L\) and \(M_k\)  can be divided into two parts. Firstly, we evolve the neurons individually by the relaxation and the external currents:
\beq \label{ODEscheme1}
\tilde{V}_{j,k}=V_{j,k}-\kl V_{j,k}-V_L\kr\Delta t+ I_0(t_k) \Delta t +\sigma_0 B_{\Delta t}^{(j)},\, j\in\{1,2,\cdots,L\},
\eeq
where the noise terms \(B_{\Delta t}^{(j)}\sim \CN(0,\Delta t)\). 

The second part  accounts for the firing events and the corresponding synaptic currents. More specifically, the scheme \eqref{ODEscheme1} probably push some \(\tilde{V}_{j,k}\) above the firing threshold \(V_F\) and causes 
\beq 
\Gamma_0^{k}=\left\{j\in\{1,2,\cdots,L\}:\tilde{V}_{j,k}\geq V_R\right\},
\eeq
to be non-empty. The neurons in \(\Gamma_0\) discharge themselves instantaneously and exert synaptic currents on the network. 


For excitatory networks where \(J>0\), we implement the cascade mechanism that accounts for spikes triggered by synaptic currents from \(\Gamma_0^k\) with an instantaneous cascade axis. Following Section 2.3, we establish a cascade axis to order the spikes iteratively
\beq 
\Gamma_i^k=\left\{j\notin\cup_{\ell=0}^{i-1}\Gamma_{\ell}^k:\tilde{V}_{j,k}+J\left|\cup_{\ell=0}^{k-1}\Gamma_{\ell}^k\right|\geq V_F\right\},\,i\in\BN^*,
\eeq
until we reach a \(\CL\in\BN^*\) such that \(\Gamma_\CL^k=\emptyset\). Then we define the cascade set
$\Gamma^k=\cup_{\ell=0}^\CL \Gamma_{\ell}^k$ indicating the neurons involved in the MFE. The update of the voltage configuration after the MFE follows either Equation \eqref{MFErule1} or Equation \eqref{MFErule2}, indicating the refractory states mechanism is not used or used respectively. If the MFE updating rule \eqref{MFErule1} is used, every neurons including neurons in \(\Gamma^k\) are impacted by the synaptic currents involved with the MFE, i.e. 
\ben\label{ODEscheme3}
V_{j,k+1}&=&\begin{cases}
\tilde{V}_{j,k}+J|\Gamma^k|,&j\notin\Gamma^k,\\
\tilde{V}_{j,k}+J|\Gamma^k|-(V_F-V_R),&j\in\Gamma^k.\end{cases}\\
M_{k+1}&=&M_k+|\Gamma^k|.\label{ODEscheme3M}
\een
If the MFE updating rule \eqref{MFErule2} is used instead, we simply reset the voltage of the neurons in $\Gamma^k$ to $V_R$, while the rest remain unchanged, i.e.
\ben\label{ODEscheme4}
V_{j,k+1}&=&\begin{cases}
\tilde{V}_{j,k}+J|\Gamma^k|,&j\notin\Gamma^k,\\
V_R,&j\in\Gamma^k.\end{cases}\\
M_{k+1}&=&M_k+|\Gamma^k|.\label{ODEscheme4M}
\een
The inhibitory networks can be treated similarly, and we skip the detailed discussions. 

 At last, we define  the firing rate in the microscopic scheme
\beq \label{dMFR}
\tilde{N}_k = \frac{M_k-M_{k-1}}{L \Delta t}.
\eeq
Although $\tilde{N}_k$ is in fact a random variable measuring the firing frequency, we can view it as an approximation to the average firing intensity when $L$ is sufficiently large. Based on the discussion above, we summarize the one-step microscopic scheme in Algorithm \ref{alg:micro}.
\begin{algorithm}
\caption{The microscopic scheme} \label{alg:micro}
\begin{algorithmic}[1]
\State \textbf{Evolving the dynamics.} With given \(\left\{V_{j,k}\right\}_{j=1}^L\) and \(M_k\), compute \(\left\{\tilde{V}_{j,k}\right\}_{j=1}^L\) with the Euler-Maruyama scheme \eqref{ODEscheme1}.
\State \textbf{Checking the firing events.} Compute the set $\Gamma^k_0$.  If it is empty, take \(V_{j,k+1}=\tilde{V}_{j,k}\) and \(M_{k+1}=M_k\) and skip to the next iteration. If not, go to Step 3.
\State \textbf{Simulating the firing events.} Compute the cascade set \(\Gamma^k\). If the refractory state mechanism is not considered, update \(V_{j,k+1}\) and \(M_{k+1}\) according to Equation \eqref{ODEscheme3} and \eqref{ODEscheme3M}. Otherwise, update \(V_{j,k+1}\) and \(M_{k+1}\) according to Equation \eqref{ODEscheme4} and \eqref{ODEscheme4M}.
\end{algorithmic}
\end{algorithm}

\subsection{The multi-scale scheme}

\newcommand{\non}{N_{\text{on}}}
\newcommand{\noff}{N_{\text{off}}}
\newcommand{\kbk}{k_{\text{back}}}

The core of the multi-scale solver is the hybridizing method combining the macroscopic scheme and the microscopic scheme. The mean firing rate \(N(t)\) which depicts the spiking activity serves as a feasible criterion to select the more suitable scheme: we apply the more efficient macroscopic scheme when the mean firing rate of the network is low, but switch to the more reliable microscopic scheme when the firing rate is relatively high, which indicates that the networks become excessively synchronous.

We partition the temporal computation domain \([T_{\min},T_{\max}]\) into \(n_t\) sub-intervals with size \(\Delta t=\frac{T_{\max}}{n_t}\), then the grid points are \(t_k=k\Delta t\). We assume the macroscopic and the microscopic schemes adopt the same temporal grid  for simplicity. To determine the more appropriate scheme in a dynamic fashion, we choose two firing rate thresholds \(\non\) and $\noff$. Recall that $N_k$ is the firing rate in the macroscopic scheme, and $\tilde N_k$ as in \eqref{dMFR} is the firing rate in the microscopic scheme. When the macroscopic scheme is in use, we continue to use it until $N_k$ exceeds \(\non\), and we switch to the microscopic solver.  When the microscopic solver is in use, we continue to use it until $\tilde N_k$ drops below \(\noff\), and we resume the use of the macroscopic solver. 

To minimize the artificial effects of the \emph{ad hoc} threshold parameters, we also introduce a buffer treatment. We choose another parameter \(\kbk\in \mathbb N\) indicating the size of the buffer zone.  When we continue adopting the macroscopic scheme until we arrive at a time \(t=t_{k_1}\) with \(N_{k_1} > \non\), we trace back to time \(t=t_{k_1-\kbk}\) and switch to the microscopic scheme.  On the other hand, We keep using the microscopic scheme until at a time \(t=t_{k_2}\) with \(\tilde N_{k} < \noff\) for every \(k=k_2-\kbk,k_2-\kbk+1,\cdots,k_2\), and we switch back to the macroscopic solver. Note that, although the buffer treatments are different in these two directions, they help ensure that the transformation between the macroscopic variables and the  microscopic variables is performed when firing activity is low in a neighborhood of the current time.

Next, we present the switching strategies between the density function and the voltage configuration. At time \(t=t_{k_1-\kbk}\), with the density function \(\left\{p_{i,k_1-\kbk}\right\}_{i=1}^{n_v}\) and the mean firing rate \(N_{k_1-\kbk}\) previously computed, we aim to transform them into a sample set of voltages \(\left\{V_{j,k_1-\kbk}\right\}_{j=1}^L\). Since at \(t=t_{k_1-\kbk}\) the firing rate  is relatively low and the network is homogeneous, we take \(L\) i.i.d. samples from \([V_{\min},V_F)\) according to the following piecewise uniform density function constructed from the discrete density \(\left\{p_{i,k_1-\kbk}\right\}_{j=1}^{n_v}\):
\beq \label{PDEtoSDE}
\rho_{k_1-\kbk}(v)=\sum_{i=1}^{n_v-1} p_{i,k_1-\kbk}\chi_{[v_{i-1/2},v_{i+1/2}]},
\eeq
where \(\chi(\cdot)\) stands for the indicator function. 
Then the microscopic scheme is reinitialized with these \(L\) samples as the initial voltage configuration \(\left\{V_{j,k_1-\kbk}\right\}_{j=1}^L\) and the counting process is reset, i.e. \(M_{k_1-\kbk}=0\). 

On the other hand, assume we have the voltage configuration \(\left\{V_{j,k_2}\right\}_{j=1}^L\) at time \(t=t_{k_2}\), and we seek to reconstruct the macroscopic variables from these samples. We transform the \(L\) voltages into a histogram and denote the number of neurons whose voltages lie in the interval \([v_{i-1/2},v_{i+1/2}]\) by \(L_i\,(i=1,2,\cdots,n_v)\). Note that we neglect the neurons whose voltages are out of computational domain of the macroscopic scheme \([V_{\min}+\frac{\Delta v}{2},V_F-\frac{\Delta v}{2}]\). The errors caused by this simplification is acceptable provided that we take \(\Delta v\) sufficiently small. Since and \(L\gg 1\), the histogram approximates the voltage density to a great extent and the discrete density follows \(p_{i,k_2}\Delta v=\frac{L_i}{\sum_{i=1}^{n_v-1}L_i}\). Finally, we arrive at the following cell averages of the density function:
\beq \label{SDEtoPDE}
p_{i,k_2}=\frac{L_i}{\Delta v\kl\sum_{i=1}^{n_v-1}L_i\kr},\quad i=1,2,\cdots,n_v-1,
\eeq
and we enforce $p_{0,k_2}=p_{n_v,k_2}=0$. Clearly, the  normalization condition holds \(\sum_{i=1}^{n_v} p_{i,k_2} \Delta v =1.\)
The  firing rate \(N_{k_2}\) is computed via the finite difference approximation:
\beq \label{SDEtoPDE2}
 N_{k_2}=\frac{M_{k_2}-M_{k_2-1}}{L \Delta t}.
\eeq

\section{Error and complexity analysis}

Throughout this section, we are concerned with approximating a \(L\) neurons integrate-and-fire network governed by Equation \eqref{sde1}-\eqref{eq:jump1}. Since the multi-scale solver evolves both the voltage configuration and the network density, we introduce a voltage dependent observable \(F: (-\infty,V_F] \rightarrow \mathbb R\)  to measure the accuracy of the microscopic solver and the macroscopic solver together. For a neuron with voltage \(V\), the observable takes value \(F(V)\). Then the ensemble mean of the observable with respect to the voltage configuration \(\{V_{j}\}_{j=1}^L\) and the density \(p(v)\) for a network are given as follows, respectively:
\beq \label{OBSV}
\mathcal{F}\kl\left\{V_{j}\right\}_{j=1}^L\kr=\frac{1}{L}\sum_{j=1}^L F(V_j),\quad
\bar{\mathcal{F}}\kl p(v)\kr=\int_{-\infty}^{V_F}p(v)F(v)\dd v.
\eeq
Similarly, the observable mean with respect to the numerical approximation of the density at the grid points is given by :
\beq 
\bar{\mathcal{F}}\kl \left\{p_i\right\}_{i=1}^{n_v-1}\kr=\sum_{i=1}^{n_v-1}p_i\int_{v_{i-1/2}}^{v_{i+1/2}}F(v)\dd v.
\eeq
We stress that the expectations of \(\mathcal{F}\kl\left\{V_{j}\right\}_{j=1}^L\kr\) and \(\bar{\mathcal{F}}\kl \left\{p_i\right\}_{i=1}^{n_v-1}\kr\) are equal provided that the voltage configuration \(\{V_{j}\}_{j=1}^L\) are i.i.d. samples from the piecewise uniform density \(\rho(v)=\sum_{i=1}^{n_v-1} p_i \chi_{[v_{i-1/2},v_{i+1/2}]}\) as in the switching scheme (Equation \eqref{PDEtoSDE}):
\beq\label{EqEXP}
\BE_{\left\{V_{j}\right\}_{j=1}^L\sim\rho}\left[\frac{1}{L}\sum_{j=1}^L F(V_j)\right]=\BE_{V\sim\rho} F(V)=\sum_{i=1}^{n_v-1}\int_{v_{i-1/2}}^{v_{i+1/2}}F(v)p_i\dd v=\bar{\mathcal{F}}\kl \left\{p_i\right\}_{i=1}^{n_v-1}\kr.
\eeq

\newcommand{\bCF}{\bar{\CF}}

To investigate the accuracy of the multi-scale solver, we study the error analysis in the following typical scenario. We adopt both microscopic and macroscopic solver to simulate the network provided the same initial data. This means the initial voltage configuration can be viewed as $L$ samples from a density function $p_0(v)$, which is also the initial condition of the macroscopic solver. The voltage configuration \(\{V_j\}_{j=1}^L\) computed via the microscopic solver  till a given time $T$ is assumed to be exact. The discrete density \(\{p_i\}_{i=1}^{n_v}\) computed via the macroscopic solver  till time $T$ is transformed into a series of i.i.d. voltage configuration \(\{\tilde{V_j}\}_{j=1}^L\) obeying Equation \eqref{PDEtoSDE} and \(\{\tilde{V_j}\}_{j=1}^L\) serves as a numerical approximation to the exact solution \(\{V_j\}_{j=1}^L\). For both cases, we have neglected the numerical error in the time integration. 

The aim of our numerical analysis is to approximate the difference between \(\{\tilde{V_j}\}_{j=1}^L\) and \(\{V_j\}_{j=1}^L\) from the perspective of computing observables, i.e.
\beq 
\CE=\left|\CF\kl\{\tilde{V_j}\}_{j=1}^L\kr-\CF\kl\{V_j\}_{j=1}^L\kr\right|.
\eeq
We assume that the exact voltage configuration \(\left\{V_{j}\right\}_{j=1}^L\) are \(L\) independent samples from a latent density function \(p(v)\), named \textbf{the prescribed density}. Using the prescribed density, the error \(\CE\) can be divided into three parts:
\ben \label{threepart} 
\CE&\leq&\CE_1+\CE_2+\CE_3\nbr\\ 
&=&\left|\CF\kl\{V_j\}_{j=1}^L\kr-\bCF(p(v))\right|+\left|\bCF(p(v))-\bCF\kl\{p_i\}_{i=1}^{n_v}\kr\right|+\left|\bCF\kl\{p_i\}_{i=1}^{n_v-1}\kr-\CF\kl\{\tilde{V_j}\}_{j=1}^L\kr\right|.\nbr
\een
Here, the error \(\CE_2\) accounts for the error of solving the Fokker-Planck equation \eqref{eq:fk} and \cite{HLXZ2021} showed that the \(L^2\) error in approximating the density is of order \(\CO\kl(\Delta v)^2\kr\). Provided that \(F\) is bounded, we immediately have \(\CE_2=\CO\kl(\Delta v)^2\kr\). The errors \(\CE_1\) and \(\CE_3\) account for the errors between i.i.d. sampled configurations and the corresponding latent densities, which are in fact random variables. The following theorem gives an estimate of \(\CE_3\) in the mean squared error.
\bthm 
Suppose that the voltage configuration \(\{\tilde{V_j}\}_{j=1}^L\) are i.i.d. samples that obey the the piecewise uniform density \(\rho(v)=\sum_{i=1}^{n_v-1} p_i \chi_{[v_{i-1/2},v_{i+1/2}]}\) as in the switching scheme Equation \eqref{PDEtoSDE}. Assume that \(F\in C_b(-\infty,V_F]\). Then the error of observables between the voltage configuration \(\{\tilde{V_j}\}_{j=1}^L\) and the density \(\{p_i\}_{i=1}^{n_v-1}\) follows
\beq 
\left[\BE_{\left\{\tilde{V_{j}}\right\}_{j=1}^L\sim\rho}\left|\bCF\kl\{p_i\}_{i=1}^{n_v-1}\kr-\CF\kl\{\tilde{V_j}\}_{j=1}^L\kr\right|^2\right]^{\frac{1}{2}}\leq \frac{C_3}{\sqrt{L}},
\eeq
where \(C_3\) is a constant relevant to \(F\) but independent of \(L\) and \(\Delta v\).
\ethm
According to Equation \eqref{EqEXP}, the expectation of \(\CF\kl\{\tilde{V_j}\}_{j=1}^L\kr\) is exactly \(\bCF\kl\{p_i\}_{i=1}^{n_v-1}\kr\) and this theorem actually copes with the variance of \(\CF\kl\{\tilde{V_j}\}_{j=1}^L\kr\).
\bprf 
 Using the independence of each \(\tilde{V_j}\), we get
\ben 
&&\BE_{\left\{\tilde{V_{j}}\right\}_{j=1}^L\sim\rho}\left|\bCF\kl\{p_i\}_{i=1}^{n_v-1}\kr-\CF\kl\{\tilde{V_j}\}_{j=1}^L\kr\right|^2\nbr\\
&=&\BE_{\left\{\tilde{V_{j}}\right\}_{j=1}^L\sim\rho}\left|\frac{1}{L}\sum_{j=1}^L\kl F(\tilde{V_{j}})-\bCF\kl\{p_i\}_{i=1}^{n_v-1}\kr\kr\right|^2\nbr\\
&=&\frac{1}{L^2}\sum_{j=1}^L\BE_{\left\{\tilde{V_{j}}\right\}_{j=1}^L\sim\rho}\left|\kl F(\tilde{V_{j}})-\bCF\kl\{p_i\}_{i=1}^{n_v-1}\kr\kr^2\right|\nbr\\
&=&\frac{1}{L}\BE_{V\sim\rho}\left|F(V)-\bCF\kl\{p_i\}_{i=1}^{n_v-1}\kr\right|^2,
\een
where we use the fact \(\BE_{V_j\sim\rho}F(V_j)=\bCF\kl\{p_i\}_{i=1}^{n_v-1}\kr\). Then we estimate the secondary moment of \(F(V)\), i.e.
\beq 
\BE_{V\sim\rho}\left|F(V)-\bCF\kl\{p_i\}_{i=1}^{n_v-1}\kr\right|^2\leq \BE_{V\sim\rho}\left|F(V)\right|^2\leq M^2,
\eeq
where \(M=\sup_{v\in(-\infty,V_F]}|F(v)|\). Take \(C_3=M\), we arrive at
\beq 
\left[ \BE |\CE_3|^2 \right]^{\frac 1 2}\leq \frac{1}{\sqrt{L}}\left[\BE_{V\sim\rho}\left|F(V)-\bCF\kl\{p_i\}_{i=1}^{n_v-1}\kr\right|^2\right]^{\frac{1}{2}}\leq \frac{C_3}{\sqrt{L}}.
\eeq
\eprf

We can give a similar estimate for \(\CE_1\) as follows:
\bthm 
Suppose that the voltage configuration \(\{V_j\}_{j=1}^L\) are i.i.d. samples that obey the the prescribed density \(p(v)\). Assume that \(F\in C_b(-\infty,V_F]\). Then the error of observables between the voltage configuration \(\{V_j\}_{j=1}^L\) and the density \(p(v)\) follows
\beq 
\CE_1=\left[\BE_{\left\{V_j\right\}_{j=1}^L\sim\rho}\left|\bCF\kl p(v)\kr-\CF\kl\{V_j\}_{j=1}^L\kr\right|^2\right]^{\frac{1}{2}}\leq \frac{C_1}{\sqrt{L}},
\eeq
where \(C_1\) is a constant relevant to \(F\) but independent of \(L\) and \(\Delta v\).
\ethm
The proof of Theorem 4.2 is rather similar to that of Theorem 4.1 and is thus omitted here.
To summarize, following Equation \eqref{threepart}, the total error \(\CE\) satisfies
\beq \label{errororder}
\CE\leq\CE_1+\CE_2+\CE_3=\CO\kl\frac{1}{\sqrt{L}}\kr+\CO\kl(\Delta v)^2\kr.
\eeq
Suppose that we evolve the macroscopic solver as an approximation to the voltage configuration of a \(L-\)neurons networks, the total error \(\CE\) involves both the numerical error of the macroscopic solver \(\CE_2\) and the errors of the switching scheme \(\CE_1\) and \(\CE_3\). As shown in Equation \eqref{errororder}, if we expect \(\CE\) to be lower than a particular threshold,  a feasible choice of spatial step length \(\Delta v\) of the macroscopic solver is of order \(\CO\kl L^{-\frac{1}{4}}\kr\), which ensures \(\CO\kl\frac{1}{\sqrt{L}}\kr=\CO\kl(\Delta v)^2\kr\). In other words, solving the Fokker-Planck equation \eqref{eq:fk} with a spatial step length of order \(\CO\kl L^{-\frac{1}{4}}\kr\) suffices to provide an accurate approximation to a network of \(L\) neurons.

Next, we compare the computational complexity of the macroscopic and the microscopic scheme. Taking the spatial step length \(\Delta v=\CO\kl L^{-\frac{1}{4}}\kr\) in the macroscopic scheme, the number of spatial grids is approximately \(\CO\kl L^{\frac{1}{4}}\kr\). For one temporary iteration, the computational complexity is \(\CO\kl L^{\frac{1}{4}}\kr\) for the explicit solver and no more than \(\CO\kl L^{\frac{3}{4}}\kr\) for the semi-implicit solver which involves solving a linear system. Compared to the microscopic scheme whose complexity for one temporal step is \(\CO(L)\), the macroscopic scheme obviously outperforms the microscopic solver. Besides, the complexity of the macroscopic solver with the semi-implicit iteration can be lowered by choosing a larger temporal step length or solving the linear system with sparse matrix techniques, since the semi-implicit method is shown to have greatly improved stability than the explicit one, see \cite{HLXZ2021}. 

Finally, we remark that the numerical analysis provided in the current work is not yet thorough due to the complicated multi-scale nature of the interacting neuron systems. But the scenario we have analyzed should be able to cover a large part of the applications, and a few other cases can be analyzed in a similar fashion. For example, we can formally derive that the error in transforming from the microscopic variables to the macroscopic variables is of the order \(\CO\kl L^{-\frac{1}{2}}\kr\) when the neurons have low correlations. Although a more comprehensive numerical analysis of the proposed method is beyond the scope of the current work, we have obtained practically useful error bounds and a proper meshing strategy for the efficient implementation of this solver. To sum up, evolving the macroscopic model as an approximation to the voltage configuration obviously lowers the computational complexity, especially when the network size \(L\gg1\), and the multi-scale solver successfully hybridizes both schemes without sacrificing the accuracy. We will further verify the efficiency of the multi-scale scheme over the microscopic scheme numerically in the next section.



\section{Numerical examples and exploration}

The numerical examples  in this section can be divided into two classes: homogeneous networks where the mean-field approximation is satisfied, and synchronous networks where the Fokker-Planck equation \eqref{eq:fk} fails to describe the original network. More specifically, we first adopt the macroscopic and microscopic solvers respectively to a weakly correlated \(L\) neurons network to verify the consistency between the microscopic interacting system \eqref{sde1}-\eqref{eq:jump1} and the macroscopic Fokker-Planck equation \eqref{eq:fk} with low firing activity. The alternating applications of the macroscopic and microscopic solvers are also tested with the switching scheme proposed in Section 3.3. Secondly, we apply the multi-scale solver to more synchronous networks, where an external pulse is exerted on the network and pushes the network to fire intensively within a period of time. We validate the performance of the proposed multi-scale solver by computing these challenging examples and comparing the results with the microscopic solver.

Without special instruction, the firing potential and the resetting potential are fixed to \(V_F=2\) and \(V_R=1\) respectively. The value \(V_{\min}\) (less than \(V_R\)) in the macroscopic scheme is adjusted in the numerical experiments to fulfill that \(p(V_{\min},t)\approx 0\) and we set \(V_{\min}=-4\). Besides, we always choose the following Gaussian initial condition for the numerical tests:
\begin{equation}\label{pg}
p_{G}(v)=\frac{1}{\sqrt{2\pi}\sigma_0}\exp\kl-\frac{(v-v_G)^2}{2\sigma_G^2}\kr,
\end{equation}
where \(v_G\) and \(\sigma_G\) are two given parameters. For the macroscopic scheme, we directly take \(p_{i,0}=p_G(v_i)\) as the initial data. For the microscopic scheme, we take \(L\) i.i.d. samples \(\left\{V_{j,0}\right\}_{j=1}^L\) from the probability density \(p_G(v)\) as the initial data. In the rest of this subsection, the mean and variance of the initial data are chosen as  \(v_G=-1\) and  \(\sigma_G^2=0.5\) respectively.

\subsection{Numerical equivalence of the microscopic/macroscopic scheme}

In this section, we are concerned with homogeneous networks with relatively low firing rates and aim to verify that the macroscopic scheme indeed provides sound approximation to the microscopic neuron networks governed by \eqref{sde1}-\eqref{eq:jump1}. The network sizes are chosen among \(L=5^4,10^4,20^4\). To avoid intensively synchronous spikes, We consider the synaptic kick \(J=-\frac{1}{L}\) for inhibitory networks and  \(J=\frac{1}{L}\) for moderately excitatory network, with the corresponding connectivity parameters fixed to \(b=-1\) and \(b=1\). The noise coefficient in Equation \eqref{sde1} is fixed to \(\sigma_0^2=\frac{1}{2}\), hence the diffusion parameter in the Fokker-Planck equation \eqref{eq:fk} follows \(a=1\). The simulation is up to time \(T=3\).

In this subsection, without special instructions, we adopt the cascade mechanism Equation \eqref{ODEscheme3} (without refractory states) when using the microscopic solver. Actually, we find that the numerical results with the cascade mechanism \eqref{ODEscheme4} (involving refractory states) are practically the same in this scenario, and are thus omitted. In other words, the numerical results strongly suggest that when the homogeneous assumptions are satisfied, whether a neuron enters a transient refractory state after a spike has little influence on the network dynamics.

To validate the accuracy of the macroscopic scheme and the switching mechanisms, we are concerned with the following four tests that evolve the neuron networks:  \vspace{-6pt}
\begin{enumerate}
\item[\textbf{\quad Test 1.}] Microscopic scheme during \([0,T]\). \vspace{-6pt}
\item[\textbf{\quad Test 2.}] Macroscopic scheme during \([0,T]\).  \vspace{-6pt}
\item[\textbf{\quad Test 3.}] Microscopic scheme during \([0,\frac{T}{2}]\) and macroscopic scheme during \([\frac{T}{2},T]\). The schemes are integrated with the switching scheme \eqref{SDEtoPDE}.  \vspace{-6pt}
\item[\textbf{\quad Test 4.}] Macroscopic scheme during \([0,\frac{T}{2}]\) and microscopic scheme during \([\frac{T}{2},T]\). The schemes are integrated with the switching scheme \eqref{PDEtoSDE}.  \vspace{-6pt}
\end{enumerate}

Among the four tests, the temporal step lengths of the macroscopic solvers and microscopic solvers involved are all fixed to \(\Delta t=5\times 10^{-5}\), small enough in order to make the time discretization error minimal. When simulating a \(L\) neurons system, the spatial step length of the macroscopic solvers involved is chosen as \(\Delta v=L^{-\frac{1}{4}}\) according to the analysis in Section 4. 

We record the density function or voltage configurations at time \(T\) computed via each Test \(n\) as \(\bigl\{p^{(n)}_i\bigr\}_{i=1}^{n_v}\,(n=2,3)\) and \(\bigl\{V_j^{(n)}\bigr\}_{j=1}^{L}\,(n=1,4)\). The voltage configurations \(\bigl\{V_j^{(n)}\bigr\}_{j=1}^{L}\,(n=1,4)\) computed via Test 1 and Test 4 are transformed into density functions \(\bigl\{p^{(n)}_i\bigr\}_{i=1}^{n_v}\,(n=1,4)\) with the switching scheme \eqref{SDEtoPDE} for sake of comparison with the density functions computed via the other two tests. Moreover, we record the firing rate at time \(t=k\Delta t\) via each Test \(n\) as \(\bigl\{N_k^{(n)}\bigr\}_{k=1}^{n_t}\). In the figures of this subsection, the firing rate of the microscopic solver is plotted per  \(k_{\textrm{rec}}\) temporal steps.



The following experiments compare the numerical results via Test 2-4 with results via Test 1, which merely rely on the microscopic model and does not adopt the coarse-grained approximation. Figure 1 and Figure 2 compare the densities and spike times computed via Test 1 and 2 for networks with different sizes \(L\) and connectivity parameters \(b\). Results in Figure 1 and Figure 2 indicate that the coarse-grained macroscopic PDE \eqref{eq:fk} gives an  accurate approximation to the microscopic network, and the accuracy is improved as \(L\) increases, which manifests the validity of mean-field approximation as \(L\to\infty\) for homogeneous networks.



\begin{figure}
\centering
\includegraphics[width=5.5cm,height=3.5cm]{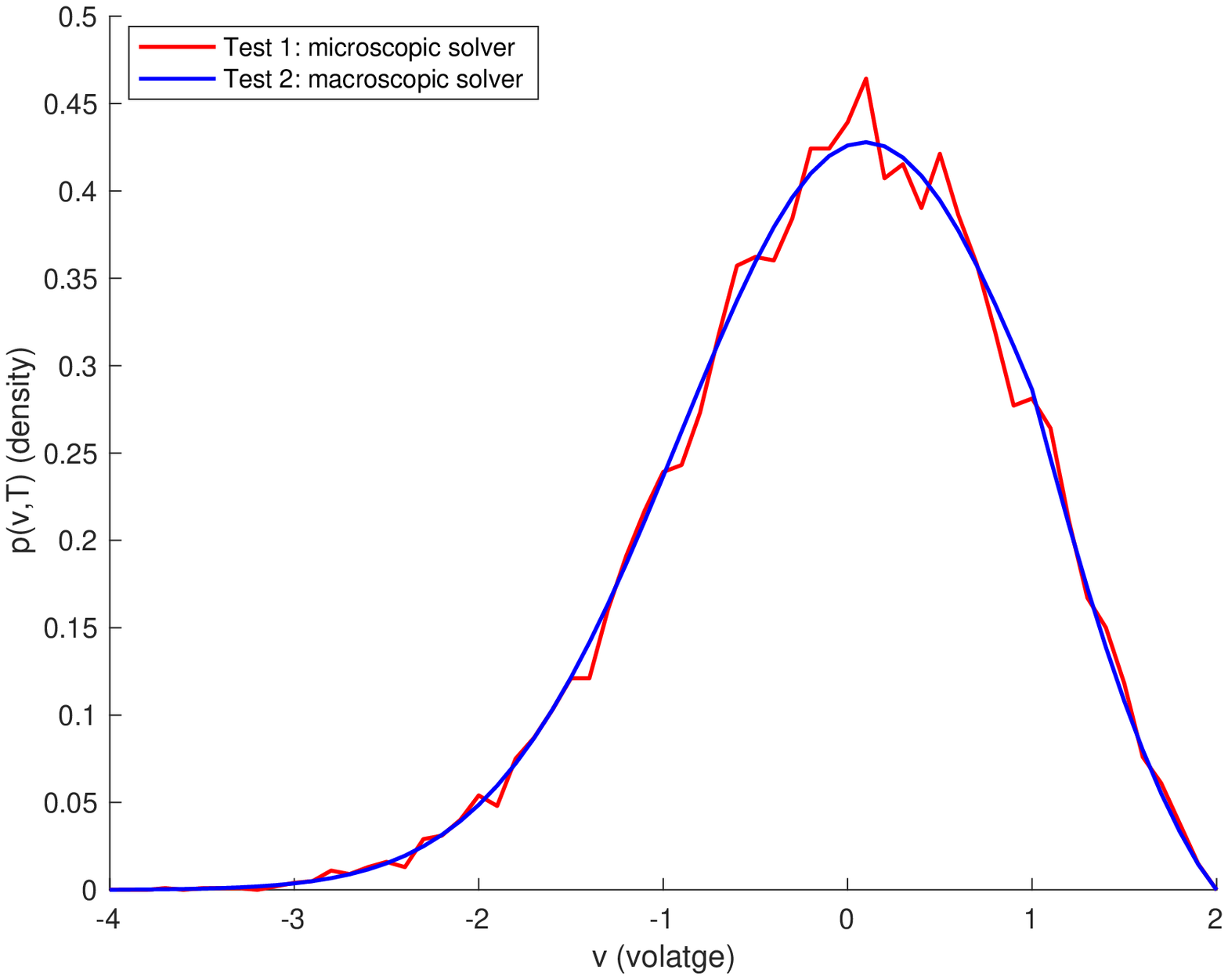}
\includegraphics[width=5.5cm,height=3.5cm]{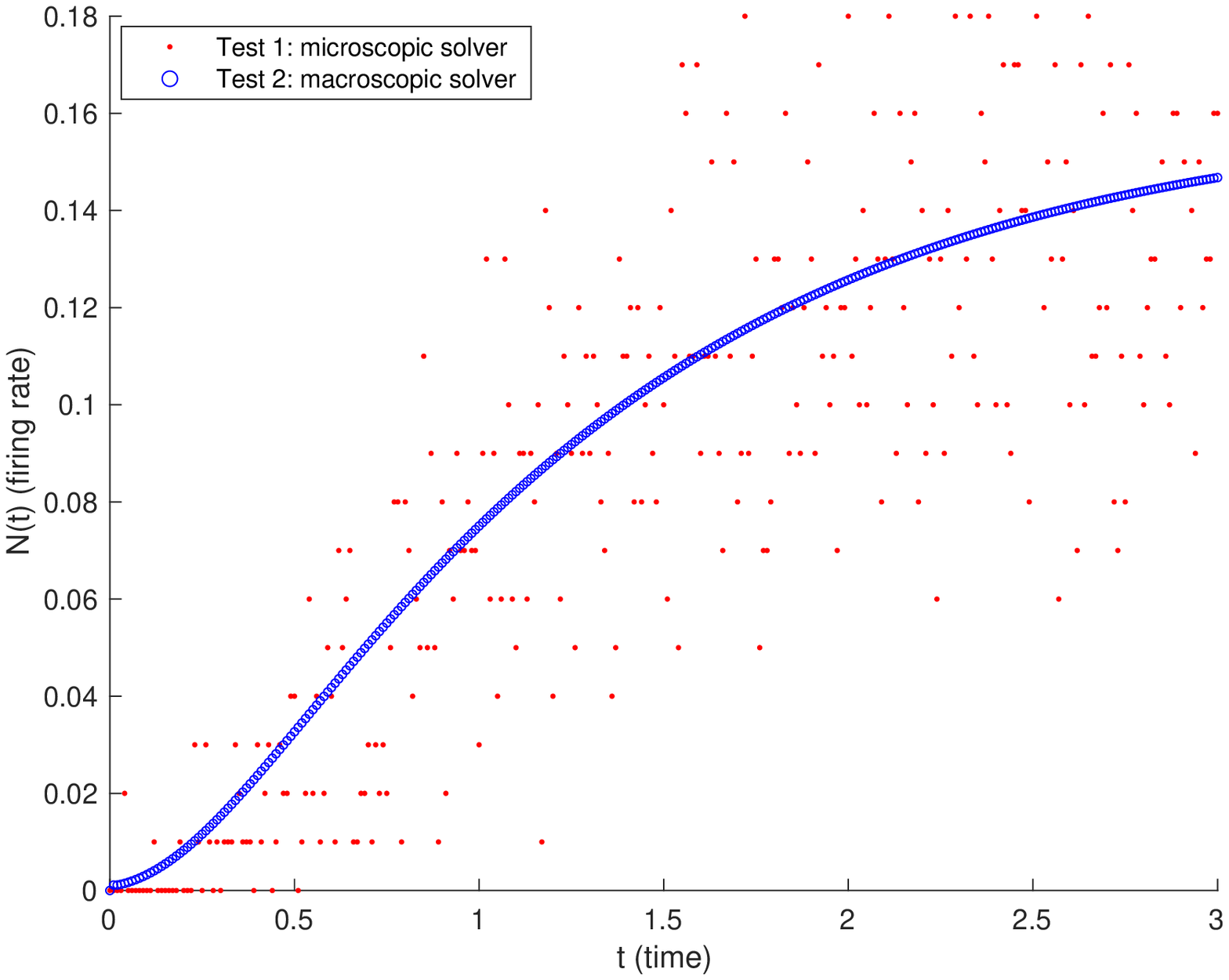}
\includegraphics[width=5.5cm,height=3.5cm]{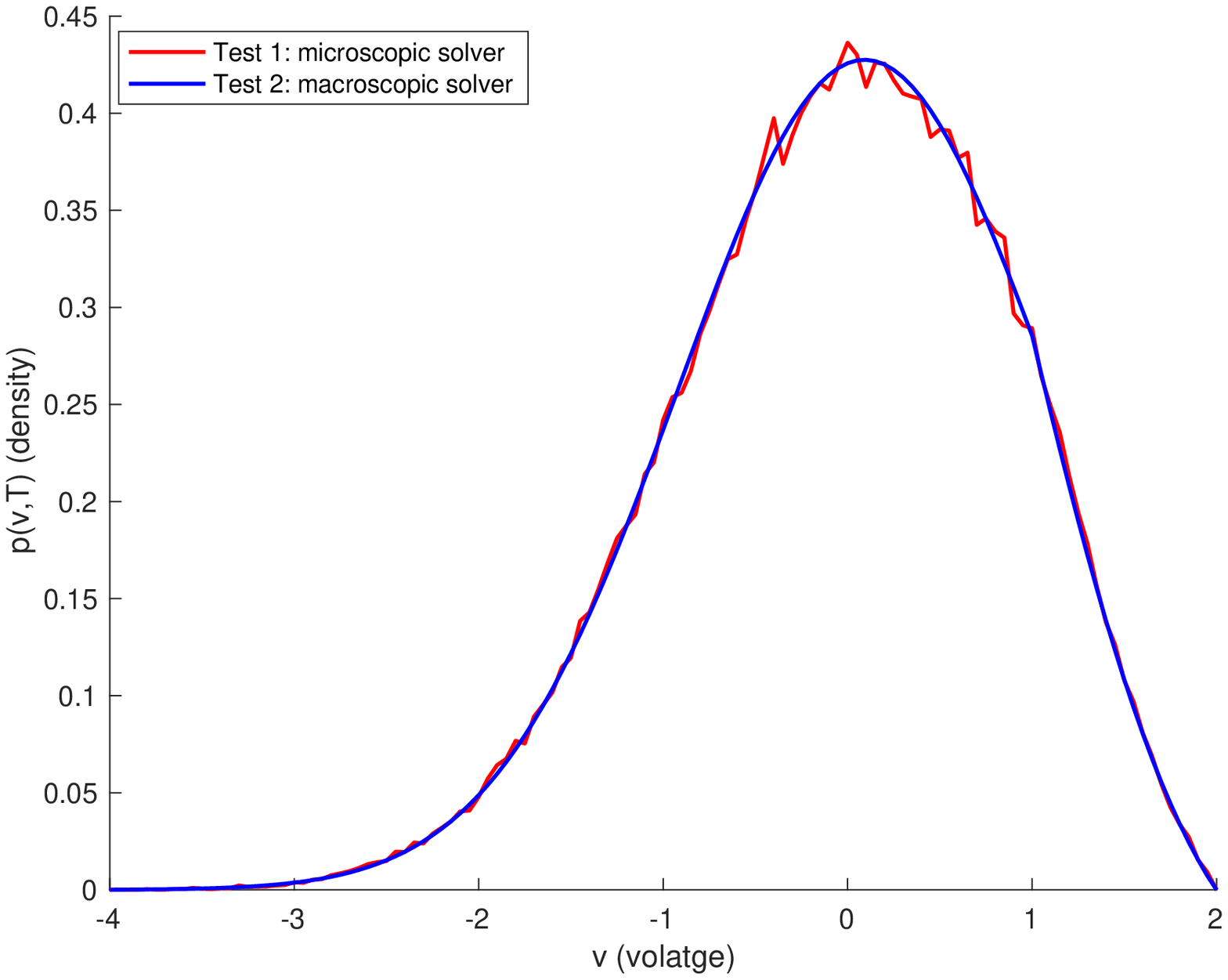}
\includegraphics[width=5.5cm,height=3.5cm]{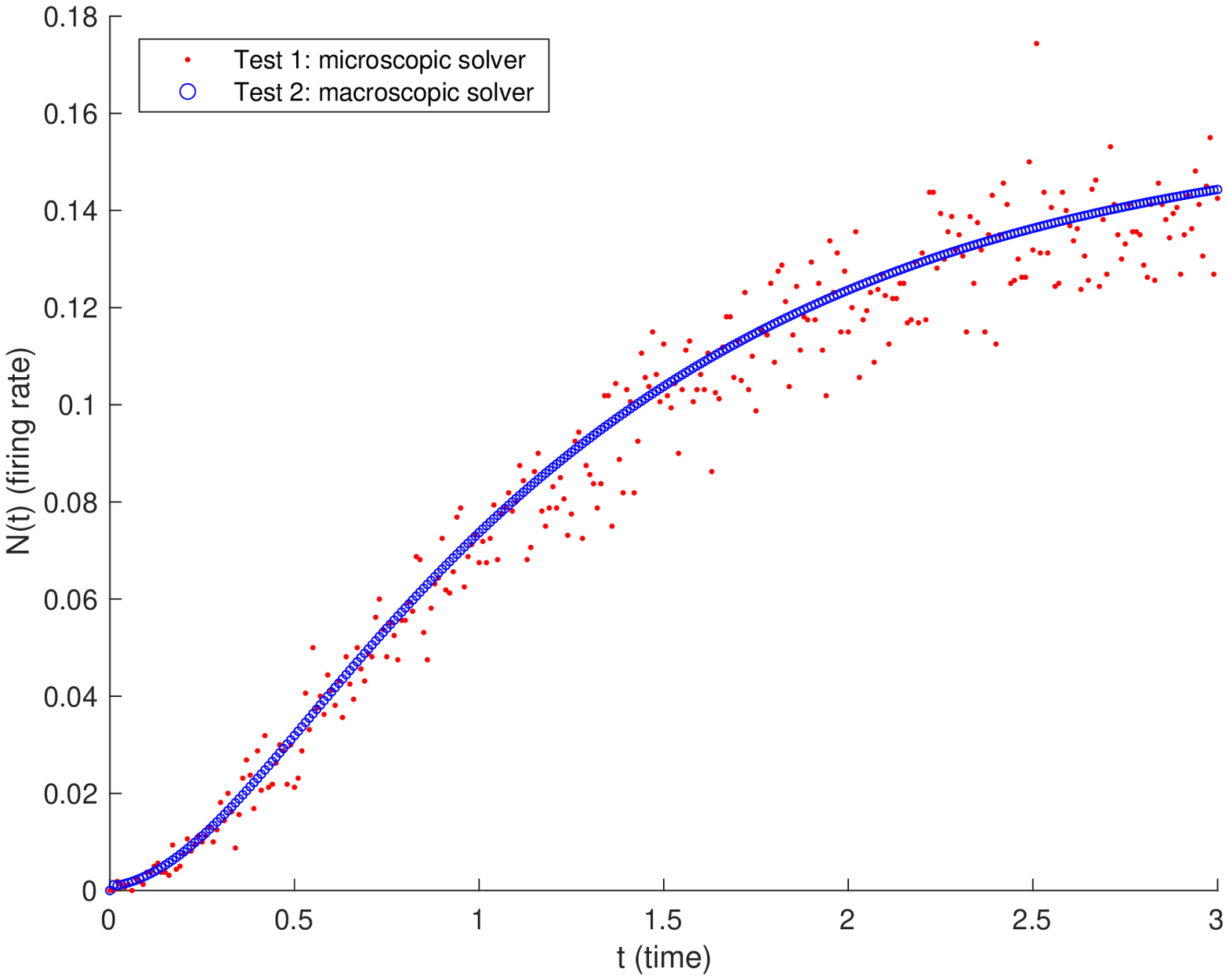}
\caption{These figures compare the network densities at time \(t=3\) and the time evolution of firing rates computed via microscopic solvers (Test 1) with results computed via macroscopic solver (Test 2). The connectivity parameter is fixed as \(b=1\). The parameter \(k_{\textrm{rec}}=100\). Upper row: network size \(L=10^4\). Bottom row: network size \(L=20^4\). Left column: network densities. Right column: firing rates.}
\label{homob1}
\end{figure}


\begin{figure}
\centering
\includegraphics[width=5.5cm,height=3.5cm]{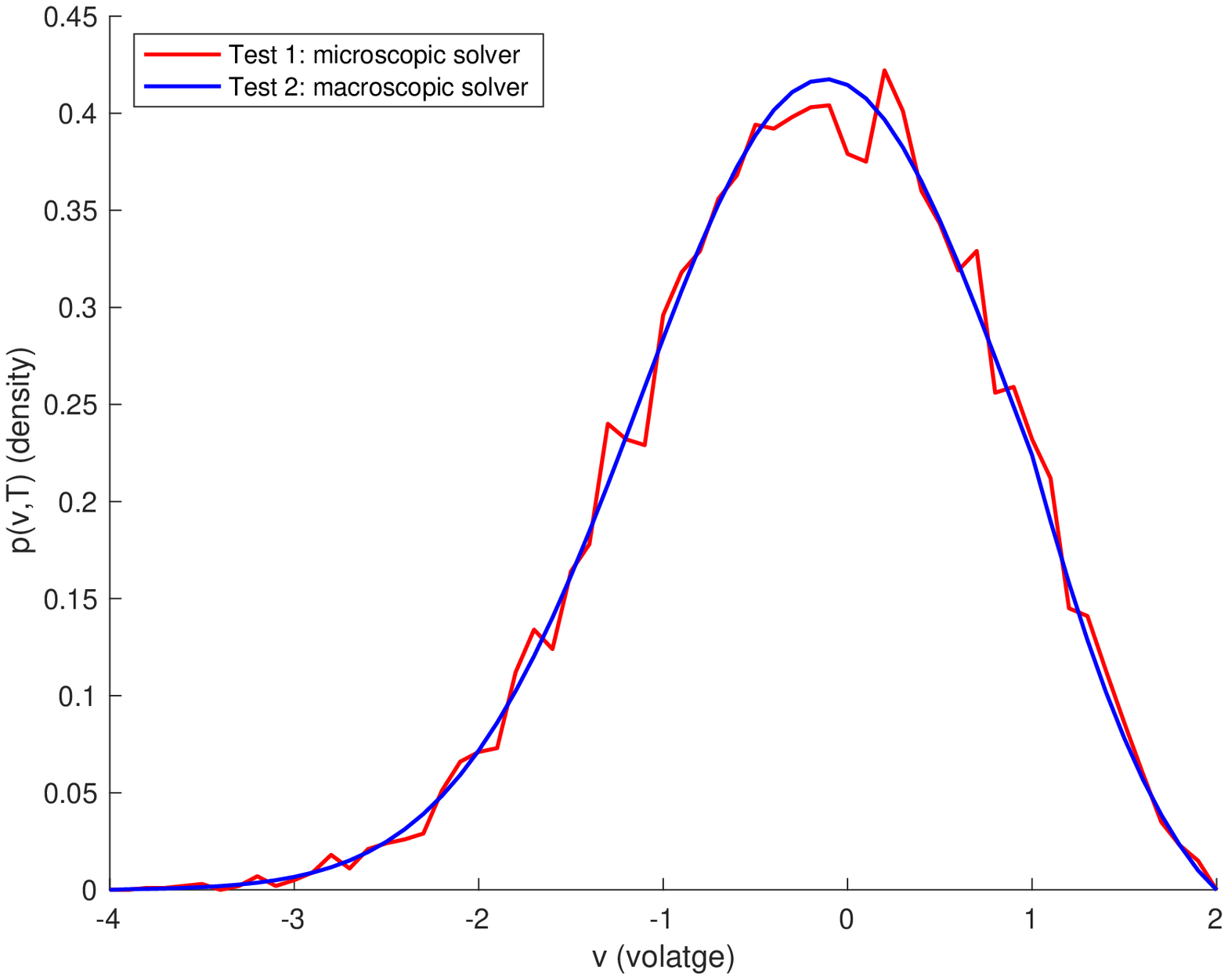}
\includegraphics[width=5.5cm,height=3.5cm]{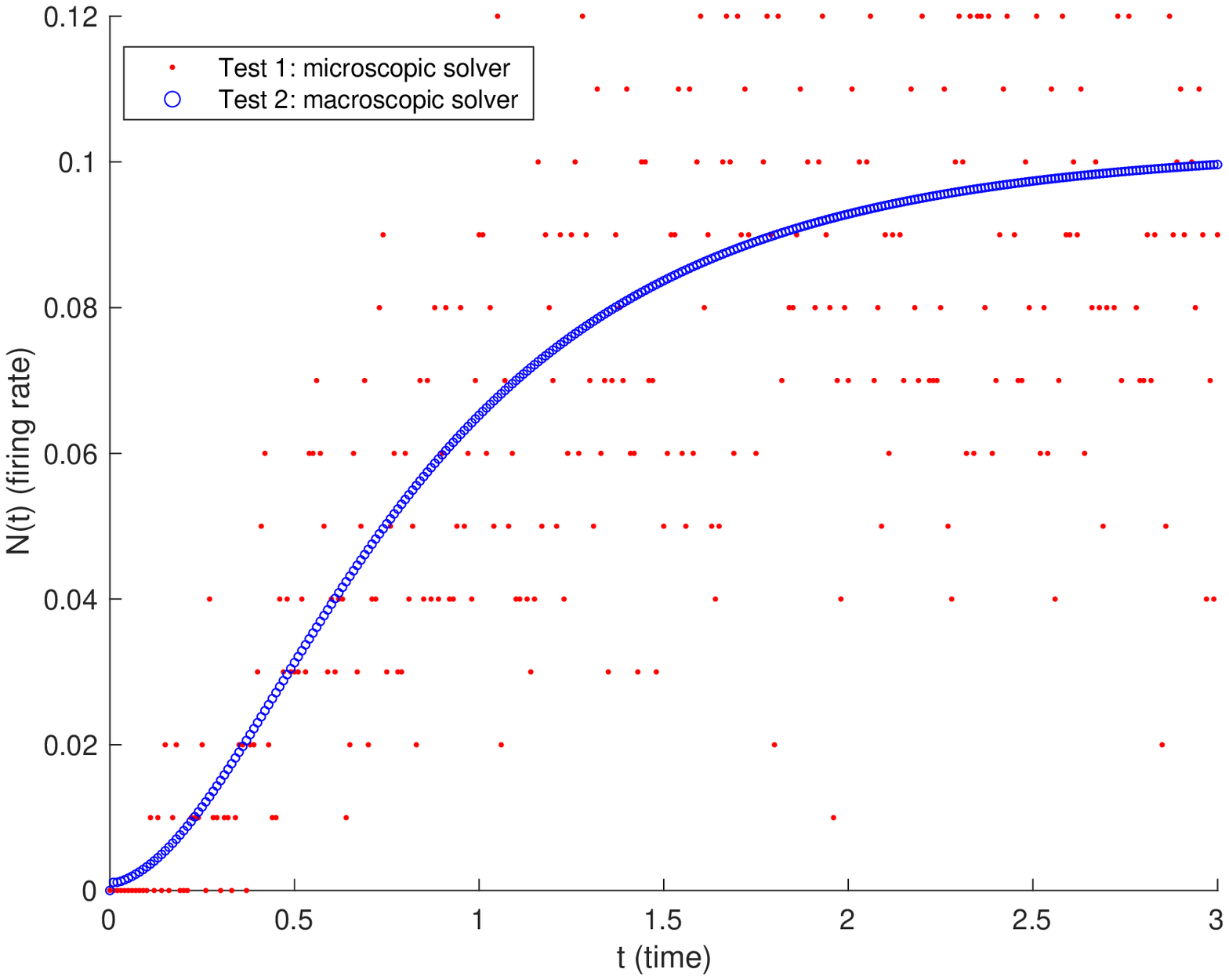}
\includegraphics[width=5.5cm,height=3.5cm]{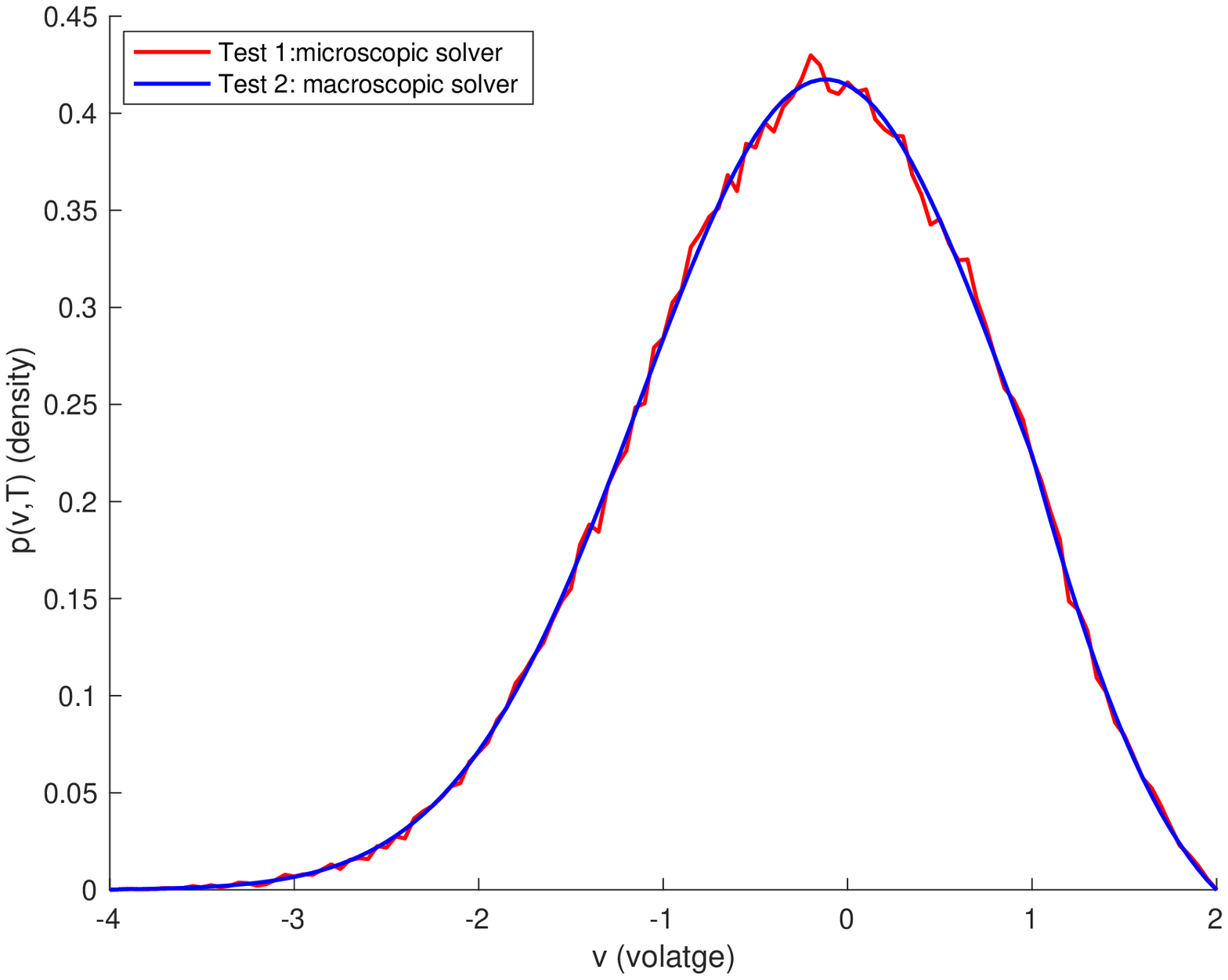}
\includegraphics[width=5.5cm,height=3.5cm]{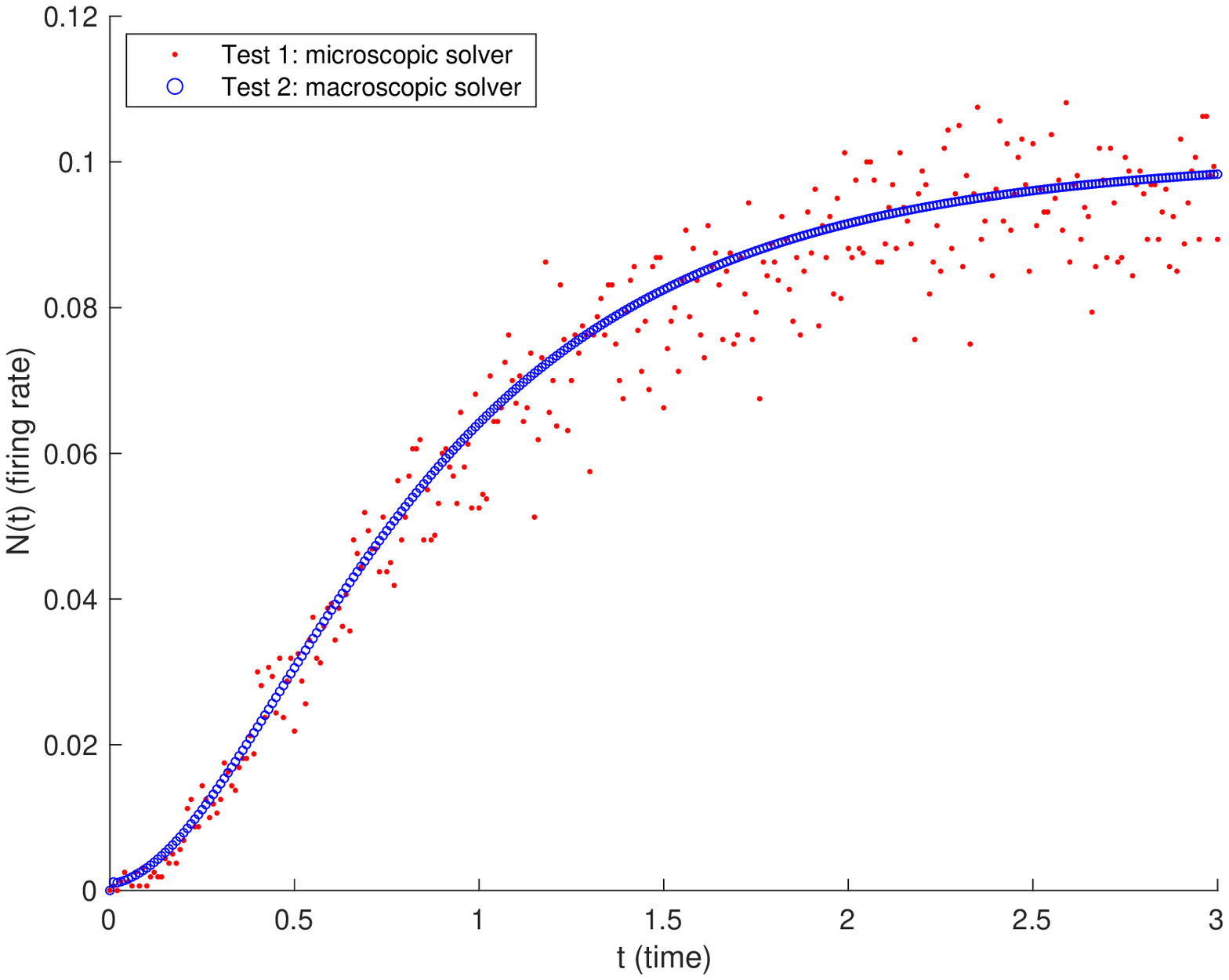}
\caption{These figures compare the network densities at time \(t=3\) and the time evolution of firing rates computed via microscopic solvers (Test 1) with results computed via macroscopic solver (Test 2). The connectivity parameter is fixed as \(b=-1\). The parameter \(k_{\textrm{rec}}=100\). Upper row: network size \(L=10^4\).  Bottom row: network size \(L=20^4\). Left column: network densities. Right column: firing rates.}
\label{homob-1}
\end{figure}



Since the microscopic solver is stochastic, we further compute the expectation of biases between network densities via different tests in order to quantize the former results. We repeat each test \(n\,(n=1,2,3,4)\) for 50 times and record the density \(\bigl\{p^{(n),l}_i\bigr\}_{i=1}^{n_v}\,(n=2,3)\) or voltage configuration \(\bigl\{V_j^{(n)}\bigr\}_{j=1}^{L}\,(n=1,4)\) for \(l=1,2,\cdots,50\). Consider three observables
\ben
F_1(v) = v, \quad 
F_2(v)=(v-\bar{v})^2, \quad 
F_3(v)=(v-\bar{v})^3.
\een
where \(\bar{v}\) denotes the the average voltage of a network, i.e. \(\bar{v}=\frac{1}{L}\sum_{j=1}^LV_j\). For density \(\bigl\{p^{(n),l}_i\bigr\}_{i=1}^{n_v}\,(n=2,3)\) or voltage configuration \(\bigl\{V_j^{(n),l}\bigr\}_{j=1}^{L}\,(n=1,4)\), the observables are computed via Equation \eqref{OBSV} and \eqref{EqEXP} respectively and recorded as \(F_1^{(n),l},F_2^{(n),l},F_3^{(n),l}\) for observable \(F_1,F_2,F_3\) respectively \((n=1,2,3,4,\, l=1,2,\cdots,50)\). Then the average biases between microscopic scheme (Test 1) and Test \(m\,(m=2,3,4)\) follows
\beq
\label{ba1}\CE_{k,m}=\frac{1}{50}\sum_{l=1}^{50}\left|F_k^{(1),l}-F_k^{(m),l}\right|,\quad k=1,\,2,\,3.
\eeq
In Table 1, we present the biases of observables of different tests with network size \(L=5^4,10^4,20^4\) respectively. Numerical results show the biases between the microscopic scheme (Test 1) and coarse-grained schemes (Test 2,3,4) decrease as \(L\) increases. Moreover, comparisons between results via Test 3 and 4 hybridizing the macroscopic and the microscopic solver and results via Test 1 provide sound evidence for the validity of the hybrid schemes. With these results above, we arrive at the conclusion that the microscopic scheme  and the hybrid schemes are numerically equivalent for networks of large sizes when firing rates are sufficiently low (i.e. homogeneous networks).

\begin{center}
  \begin{tabular}{c|ccc}
   \toprule
   \textbf{First-order observable} & Test 1 - Test 2 & Test 1 - Test 3 & Test 1 - Test 4  \\
   \midrule
   \(L=5^4\) & 3.57e-02 & 3.68e-02 & 5.02e-02\\
   \(L=10^4\) & 8.35e-03 & 8.37e-03 & 1.10e-02\\
   \(L=20^4\) & 3.11e-03 & 2.91e-03 & 2.67e-03\\
   \midrule
   \textbf{Second-order observable} & Test 1 - Test 2 & Test 1 - Test 3 & Test 1 - Test 4   \\
   \midrule
   \(L=5^4\) & 2.07e-02 & 2.11e-02 & 3.32e-02\\
   \(L=10^4\) & 6.67e-03 & 6.55e-03 & 5.73e-03\\
   \(L=20^4\) & 3.59e-03 & 3.53e-03 & 1.53e-03\\
   \midrule
   \textbf{Third-order observable} & Test 1 - Test 2 & Test 1 - Test 3 & Test 1 - Test 4 \\
   \midrule
   \(L=5^4\) & 4.52e-02 & 4.51e-02 & 6.63e-02\\
   \(L=10^4\) & 1.12e-02 & 1.12e-02 & 1.55e-02\\
   \(L=20^4\) & 3.20e-03 & 3.23e-03 & 3.69e-03\\
   \bottomrule
  \end{tabular}
 \captionof{table}{The table compares the density biases computed via different tests as shown in Equation \eqref{ba1}. The connectivity parameter is fixed as \(b=1\). The network size \(L\) varies.}
 \label{TabGau}
 \end{center}

\subsection{Numerical tests on the multi-scale scheme}

In this subsection, we further adopt the hybrid multi-scale scheme proposed in Section 3.3 on synchronous networks when neurons fire more actively. The numerical performance of the multi-scale solver on these synchronous networks is validated through comparing them with the reliable microscopic solver. We aim to show that thI’m e multi-scale solvers produce accurate numerical solutions with less simulation time compared with the microscopic scheme. 

\subsubsection{Preparation}

In this part, we specify the choice of the network parameters. We also discuss the refractory mechanism in the microscopic solver that affects the simulations of synchronous networks.

Recall that the multi-scale solver treats networks in homogeneous and synchronous regimes with different solvers, we aim to test the multi-scale solver on an "intermittent synchronous network". That is, excitatory networks that become extremely synchronous for a transient period due to the external input, but return to the homogeneous status shortly. To identify this kind of network, we need to choose moderately large \(b\) and set the external current \(I_0(t)\) to a suitable amplitude to ensure the networks emit moderately intensive MFEs.

The following test aims to find suitable \(b\) and \(I_0(t)\). Firstly, we define a time-dependent input pulse as follows
\beq \label{IbI0}
I_0(t)=\CJ(t;J_{\max},\beta,t_p)=J_{\max}\exp\left[-\beta(t-t_p)^2\right].
\eeq 
Here, \(J_{\max}>0\) and \(t_p>0\) denote the amplitude and input time of the pulse, while \(\beta\gg 1\) depicts the concentration of the pulse. 


Fixing $\beta=100$ and $t_p=0.2$, for each \(b\), we simulate the networks different input amplitude \(J_{\max}\) and find the minimum \(J_{\max}\) that makes the firing rate exceed a threshold \(N_0=15\) during a fixed finite time interval \([0,1]\). The results are shown in Figure \ref{Ib} (Left). We emphasize that the threshold \(N_0=15\) is simply a moderately large number and the macroscopic solver might not agree with the microscopic solver when $N>N_0$. When adopting the multi-scale solver, the threshold is chosen \(\non<N_0=15\) to ensure the validity of the macroscopic solver. 


\begin{figure}
\centering
\includegraphics[width=6cm,height=4cm]{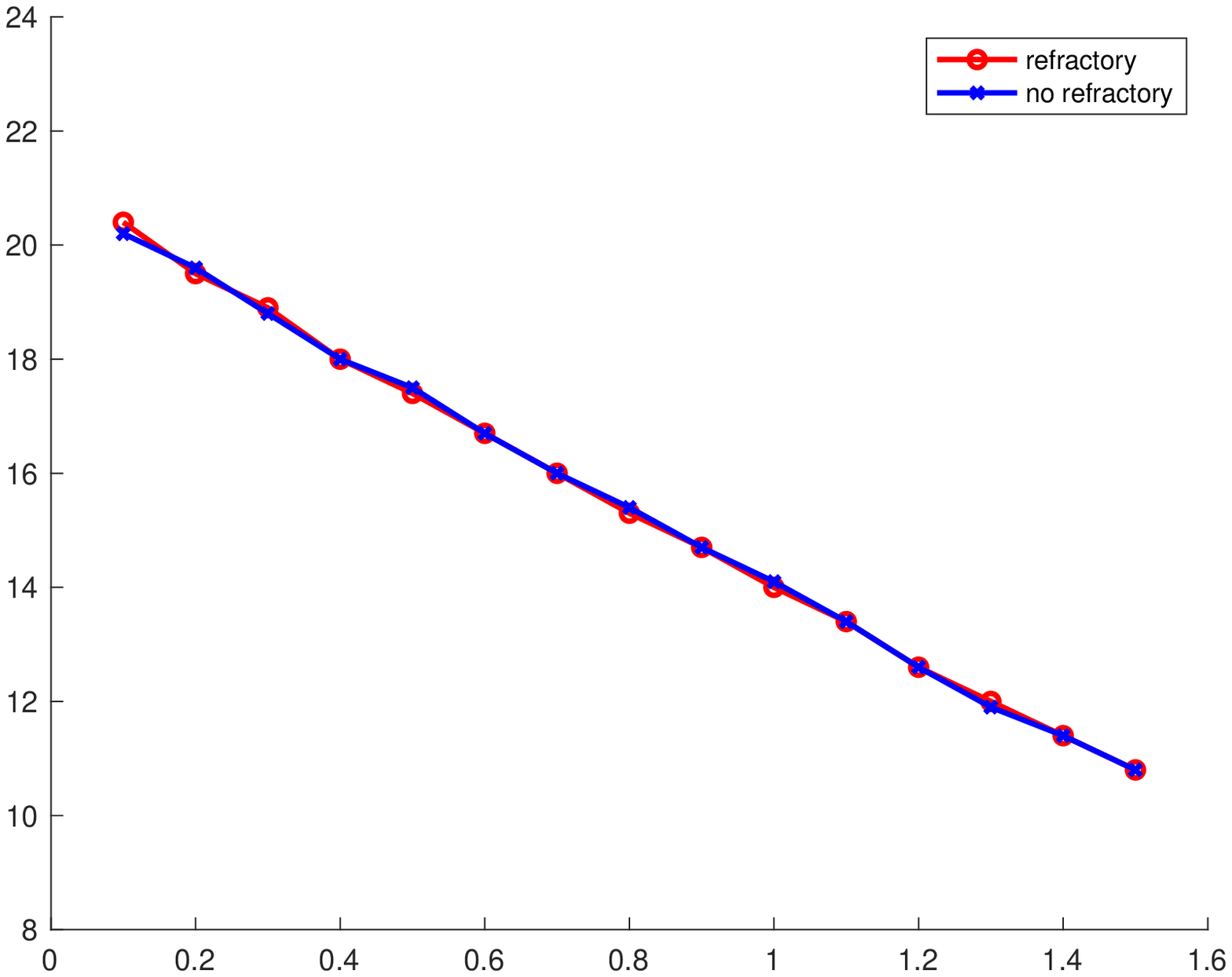}
\includegraphics[width=6cm,height=4cm]{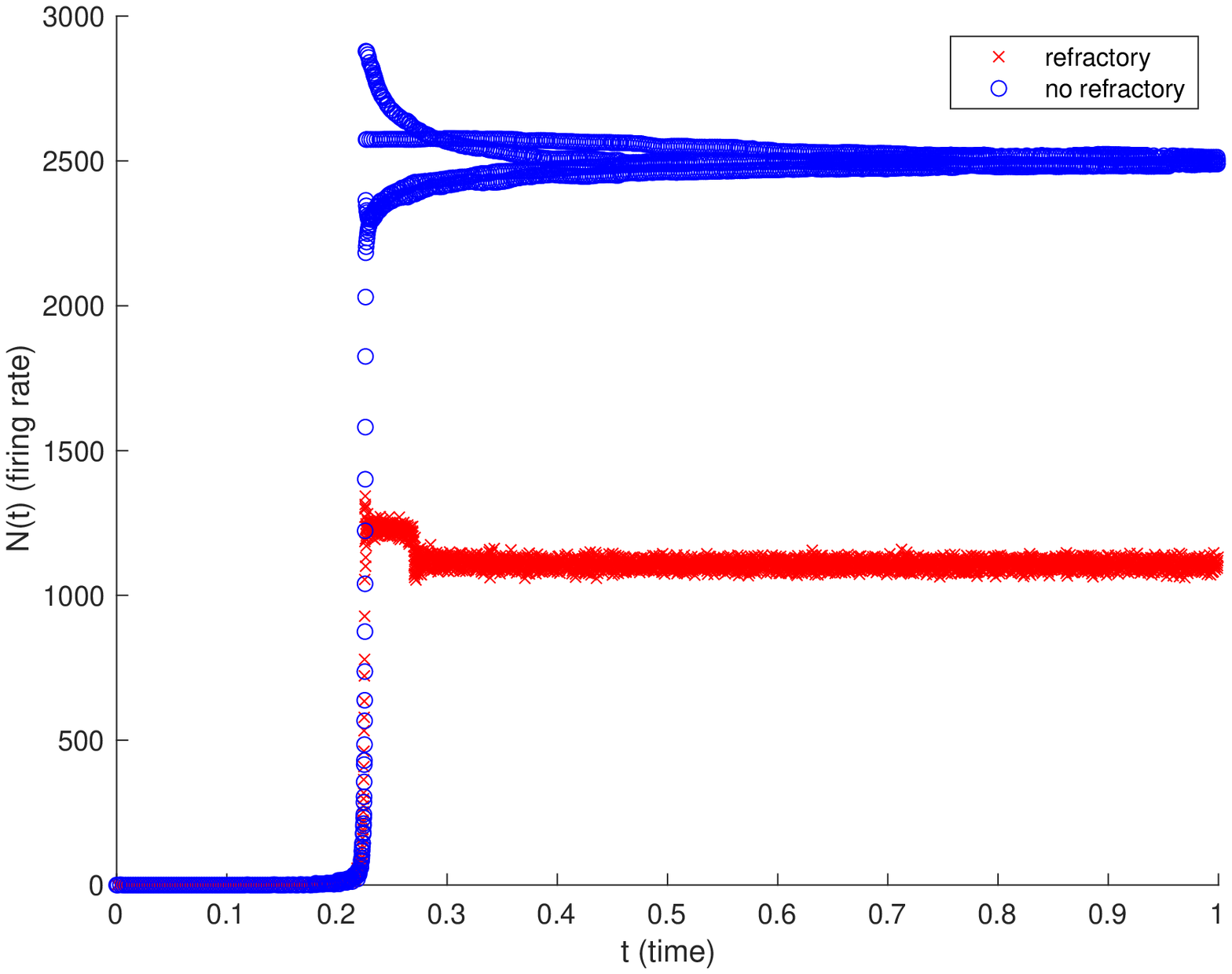}
\caption{The left figure shows the minimum amplitude \(J_{\max}\) that pushes the firing rate exceed a threshold \(N_0=15\) in a finite time interval \([0,1]\) for every fixed \(b\). The right figure show the time evolution of the firing rates of networks with \((b,J_{\max})=(1.2,20)\) computed via the macroscopic solvers with and without the refractory states mechanism. The tests are implemented on networks with sizes \(L=10^4\).}
\label{Ib}
\end{figure}


We highlight that the tests involved with Figure \ref{Ib} (Left) are compared with the microscopic scheme with and without the refractory mechanism, i.e. microscopic schemes with voltage updating rule \eqref{ODEscheme3} and \eqref{ODEscheme4} respectively. In this scenario, both schemes still yield practically the same results, implying that when $N>N_0$, the role of the refractory state is still minimal. However, if we take  \((b,J_{\max})\) far above the curve in Figure \ref{Ib} (Left) such that the  networks become highly synchronous, we observe  that the refractory mechanism  does matter  to the simulation results.   We test both microscopic solvers to neuron networks with connectivity parameters \(b=1.2\) and \(J_{\max}=20\). The networks fire extremely actively and become considerably synchronous after the intensive pulse input in \(t=0.2\). The numerical results of the synchronous networks by those two reset rules \eqref{MFErule1} and \eqref{MFErule2} are shown in Figure \ref{Ib} (Right), and they differ greatly. We observe that in this case, the role of the refractory state is no longer negligible. Since the simplified updating rule \eqref{MFErule1} may fail to produce the physical solution, we choose to use the updating rule \eqref{MFErule2} in the microscopic scheme for the rest of the section.

\subsubsection{Singe pulse Experiments}

In this part we test on neuron networks with \(I_0(t)\) in the form of single pulse, i.e.
\beq\label{T5I0}
I_0(t) = \CJ(t;J_{\max},100,0.5).
\eeq
Following the idea of finding moderately synchronous networks, we take \(b=1,J_{\max}=16\) and \(b=0.5,J_{\max}=24\) as the network parameters. These parameters ensure \((b,J_{\max})\) to be on top of the slope in Figure \ref{Ib}, but are not so large such that the networks do not fire intensively. The tests are carried on time interval \([0,1]\) using the macroscopic scheme and the multi-scale scheme with temporal step length \(\tau=10^{-4}\) respectively. The network sizes and spatial step lengths are fixed as \(L=20^4\) and \(\Delta v=\frac{1}{20}\), following the condition \(\Delta v=\CO\kl L^{-\frac{1}{4}}\kr\) as shown in Section 4. The solvers adopt the refractory mechanism. The firing rate thresholds are fixed to \(\non=\noff=10\) and the buffer zone size is chosen as \(\kbk=10\).

\begin{figure}
\centering
\includegraphics[width=5.5cm,height=3.5cm]{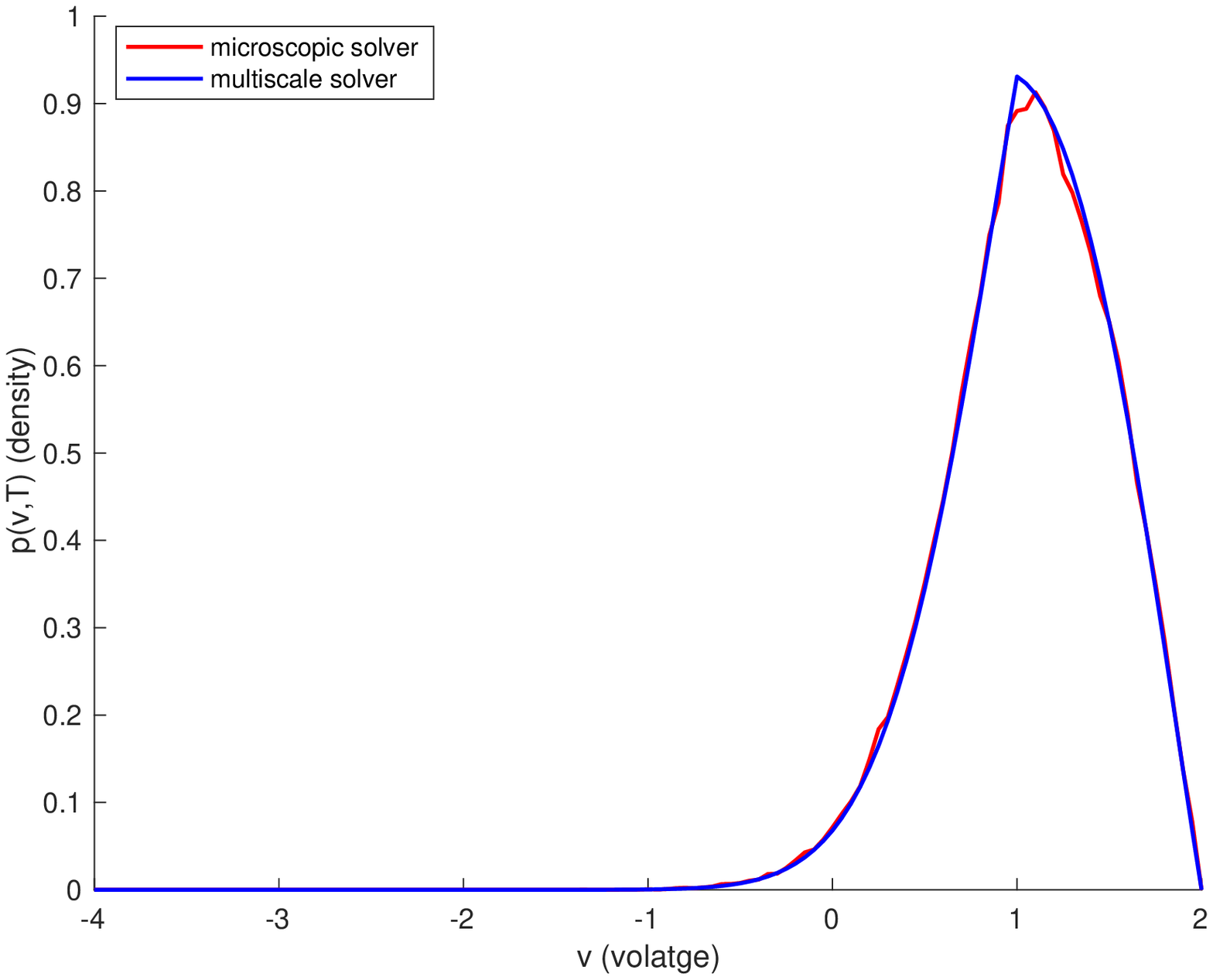}
\includegraphics[width=5.5cm,height=3.5cm]{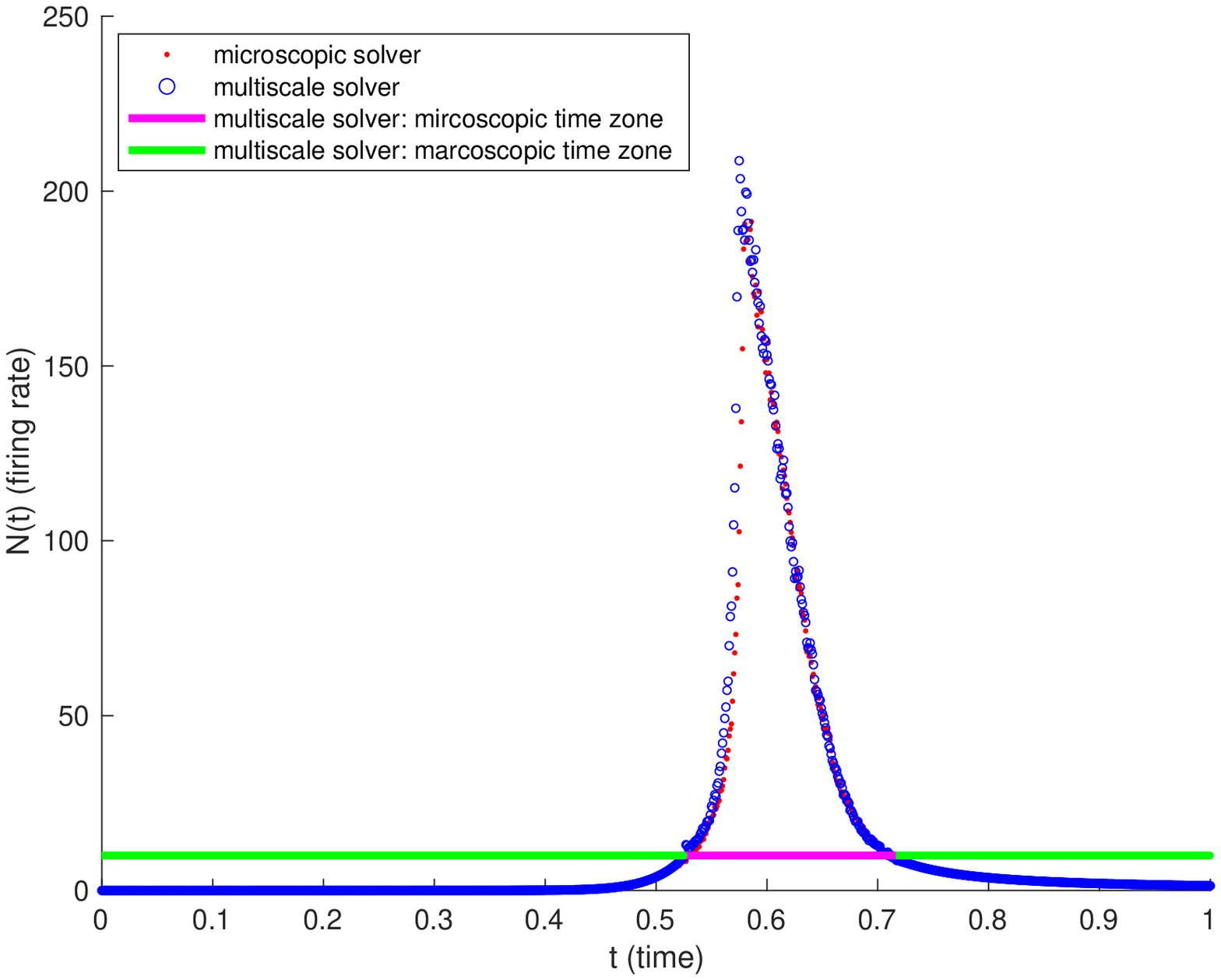}
\includegraphics[width=5.5cm,height=3.5cm]{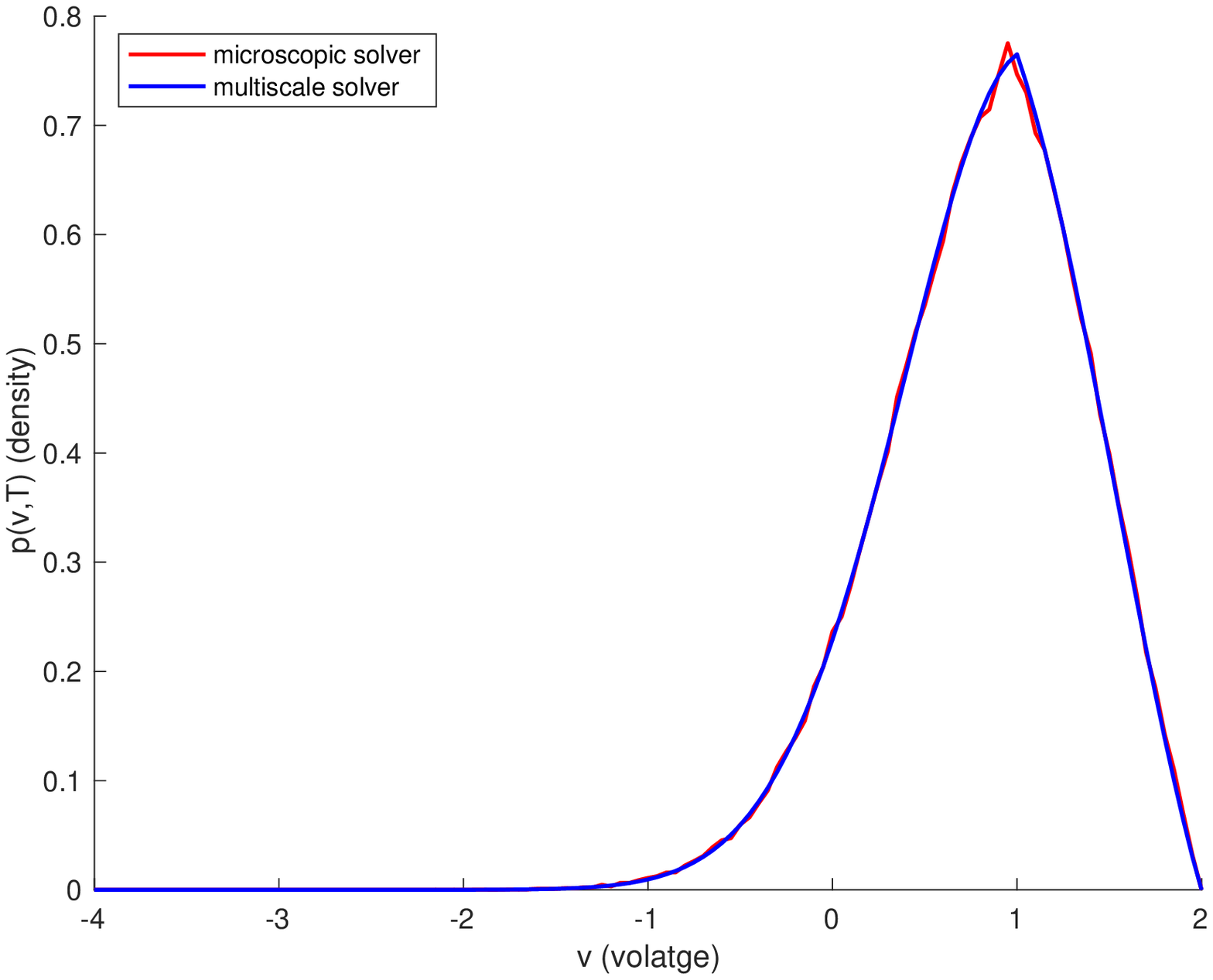}
\includegraphics[width=5.5cm,height=3.5cm]{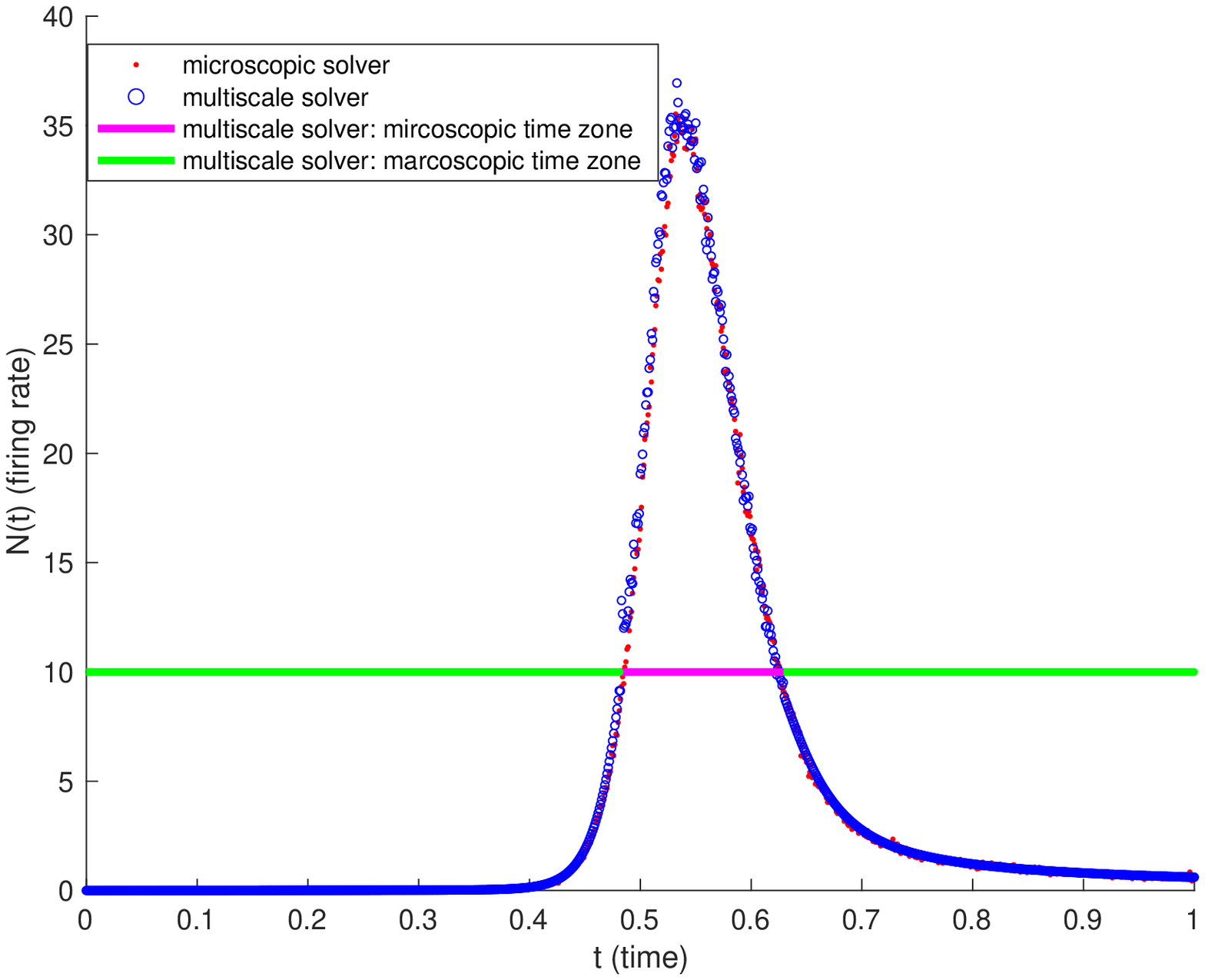}
\caption{We adopt the microscopic scheme and the multi-scale scheme on networks with input currents \eqref{T5I0}. These figures compare the densities of the voltage configurations at time \(t=1\) and the evolution of the firing rate on time interval \([0,1]\). The parameter \(k_{\textrm{rec}}=10\). Upper row: network parameters \((b,J_{\max})=(1,15)\). Bottom row: network parameters \((b,J_{\max})=(0.5,24)\). Left column: the network densities. Right column: evolution of the firing rates.}
\label{T5}
\end{figure}

\textbf{\begin{figure}
\centering
\includegraphics[width=5.5cm,height=3.5cm]{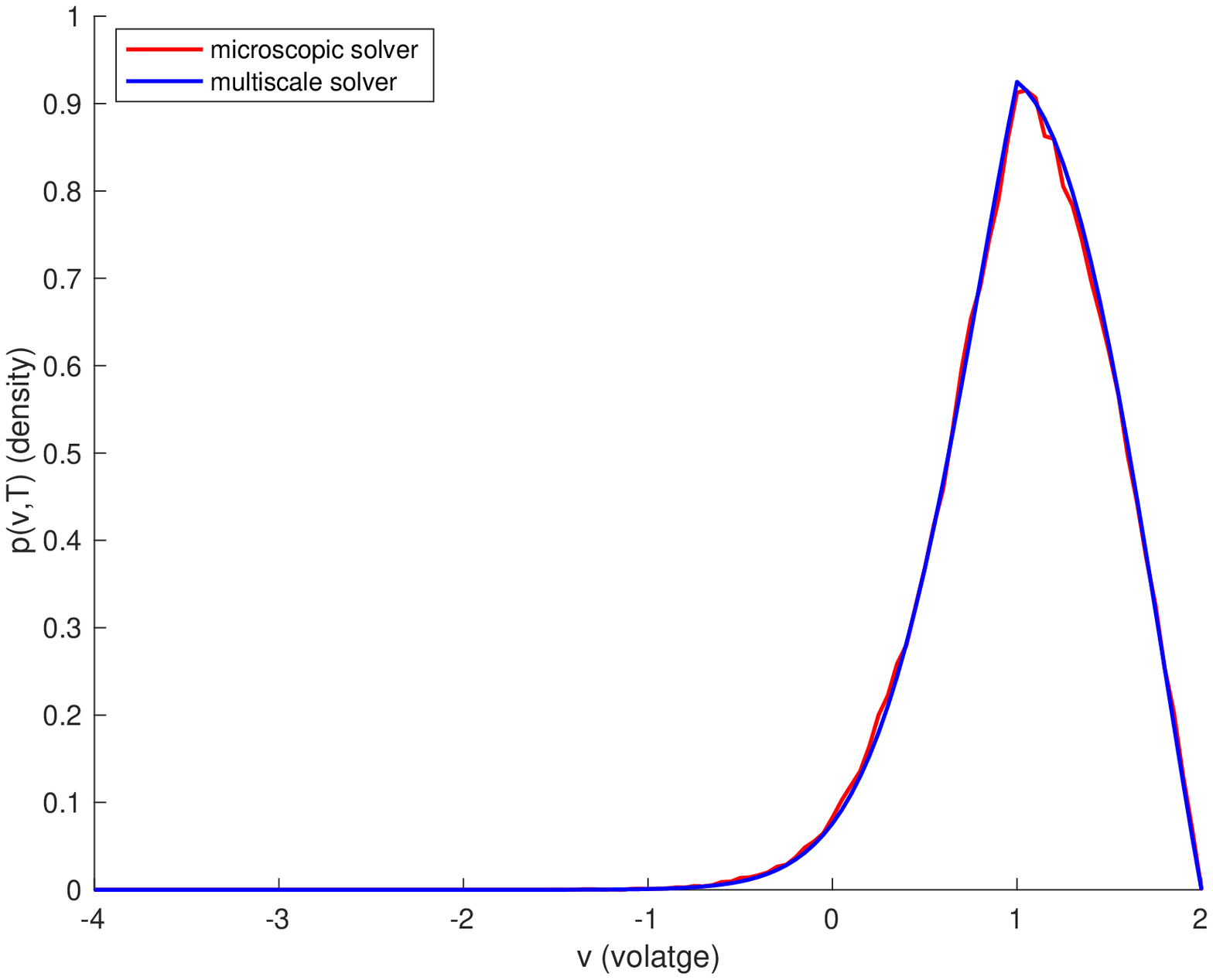}
\includegraphics[width=5.5cm,height=3.5cm]{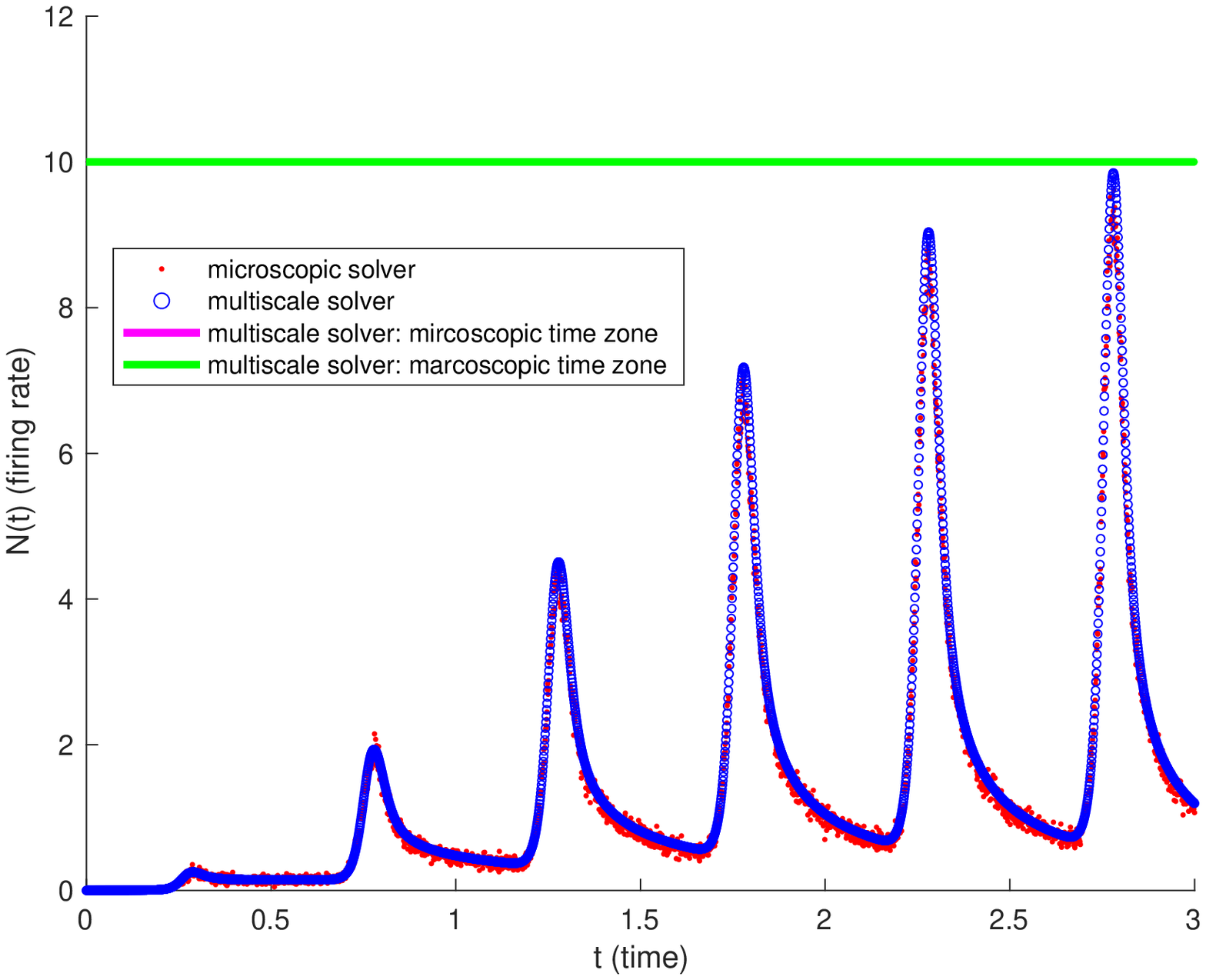}
\includegraphics[width=5.5cm,height=3.5cm]{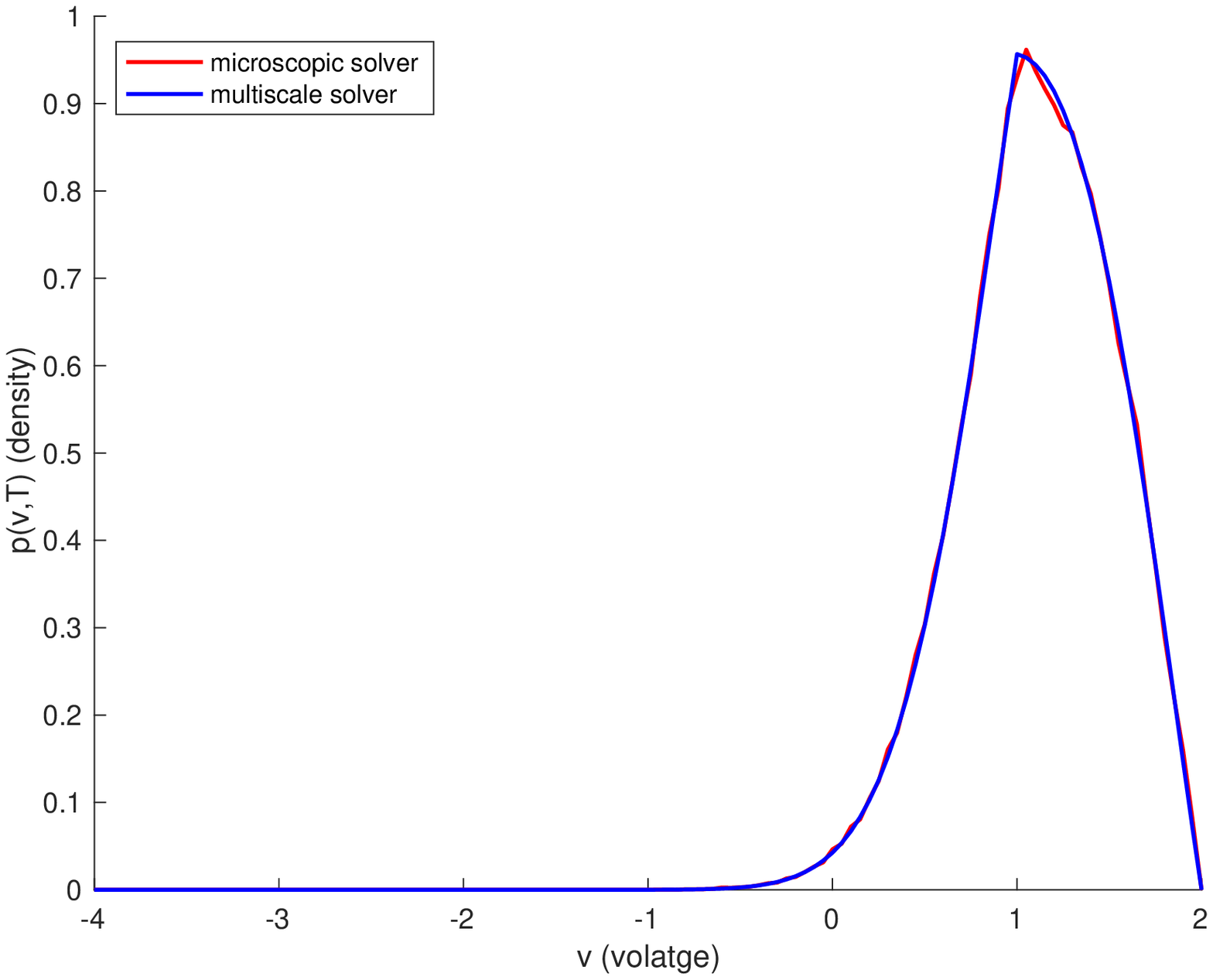}
\includegraphics[width=5.5cm,height=3.5cm]{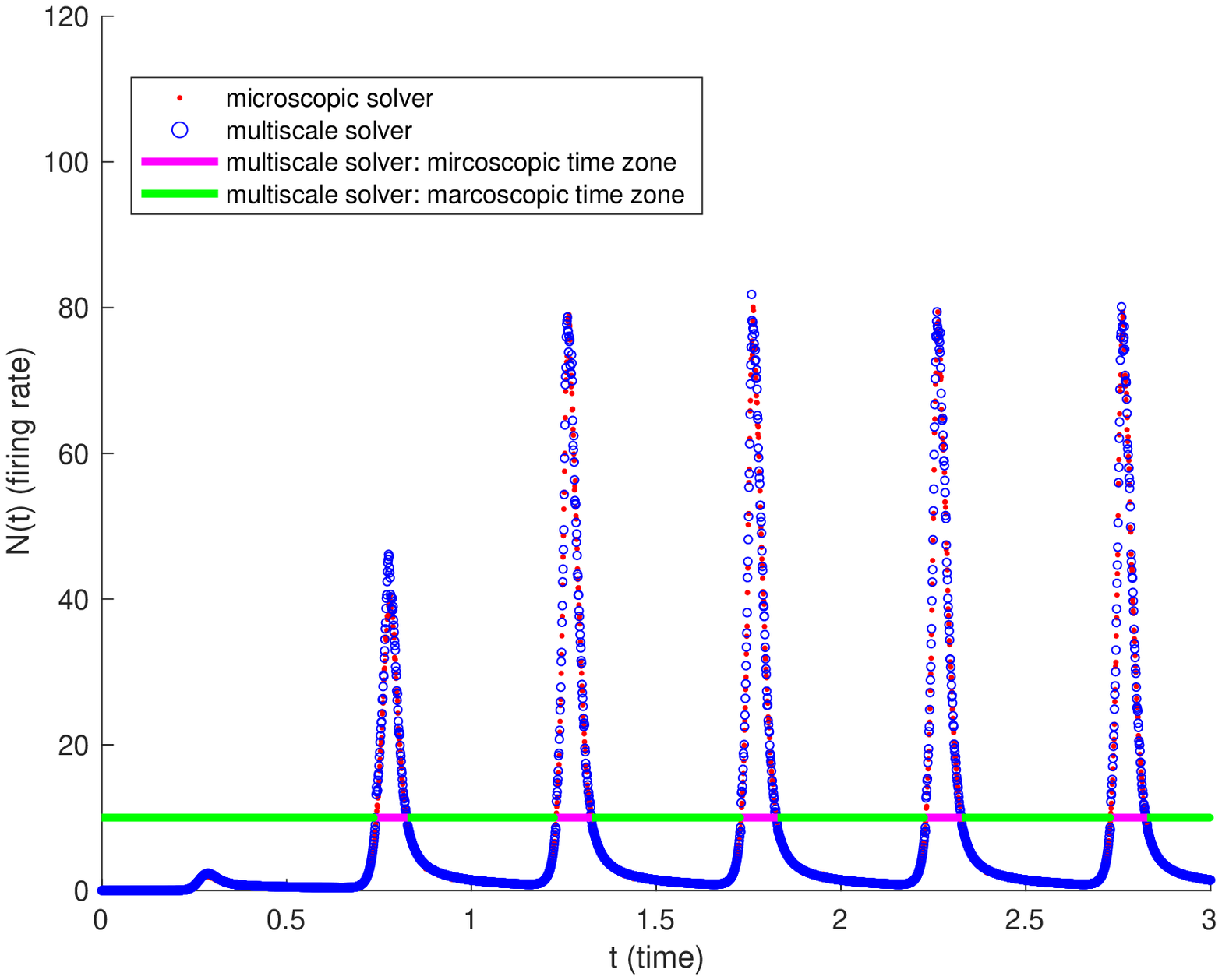}
\caption{We adopt the microscopic scheme and the multi-scale scheme on networks with input currents \eqref{T6I0}. These figures compare the densities of the voltage configurations at time \(t=3\) and the evolution of the firing rate on time interval \([0,3]\). The parameter \(k_{\textrm{rec}}=10\). Upper row: network parameters \((b,J_{\max})=(0.8,10)\). Bottom row: network parameters \((b,J_{\max})=(0.8,20)\). Left column: the network densities. Right column: evolution of the firing rates.}
\label{T6}
\end{figure}}

In Figure \ref{T5}, we show the voltage configurations and the time evolution of the firing rates computed via the microscopic solver and the multi-scale solver. The results indicate that the multi-scale solver produces a reliable approximation to the microscopic solver even if the network becomes synchronous for a transient period of time. From the figures of the firing rates, we find that the multi-scale solver merely uses the microscopic solver around time \(t=0.5\) when the pulse is input to the network, while it uses the microscopic solver at other times. Through adopting the macroscopic solver, the multi-scale solver remarkably reduces the time cost of simulating the network and we present the system times of the two schemes in Table 2.

\begin{center}
  \begin{tabular}{c|cccc}
   \toprule
   \((b,J_{\max})=(1,16)\) & \(L=10^4\) & \(L=15^4\) & \(L=20^4\) & \(L=25^4\)  \\
   \midrule
   microscopic  & 1.60 & 6.30 & 18.27  & 42.19\\
   multi-scale & 0.77 & 2.04 & 5.04 & 10.43\\
   \midrule
   \((b,J_{\max})=(0.5,24)\) & \(L=10^4\) & \(L=15^4\) & \(L=20^4\) & \(L=25^4\)  \\
   \midrule
   microscopic  & 1.75 & 6.51 & 18.16  & 41.91\\
   multi-scale & 0.78 & 1.83 & 3.84 & 7.91\\
   \bottomrule
  \end{tabular}
 \captionof{table}{The table presents the system times to simulate a \(L\) neurons network with input current \eqref{T5I0} with network parameter \((b,J_{\max})\).}
 \label{TabT5}
 \end{center}

\subsubsection{Periodic pulses Experiments}

In this part we test on neuron networks with periodic synaptic currents \(I_0(t)\) that is a linear combination of several single pulses, i.e.
\beq\label{T6I0}
I_0(t) = \sum_{k=0}^5\CJ(t;J_{\max},500,\frac{k}{2}+\frac{1}{4}).
\eeq
In the following experiments, we fixed the parameter \(b=0.8\) and choose two input amplitude \(J_{\max}=10,20\) as the network parameters to yield moderately synchronous networks. The pulse concentration in input current \eqref{T6I0} is made higher to avoid interference between each input pulse. The tests are carried out on time interval \([0,3]\) using the microscopic scheme and the multi-scale scheme with temporal step length \(\tau=10^{-4}\) respectively. The network sizes and spatial step lengths are fixed as \(L=20^4\) and \(\Delta v=\frac{1}{20}\), following the condition \(\Delta v=\CO\kl L^{-\frac{1}{4}}\kr\) as shown in Section 4. The solvers adopt the refractory mechanism. The firing rate thresholds are fixed to \(\non=\noff=10\) and the buffer zone size is chosen as \(\kbk=10\).

In Figure \ref{T6}, we show the voltage configurations and the time evolution of the firing rates computed via the microscopic solver and the multi-scale solver. When using the multi-scale solver, the firing rate in the \((b,J_{\max})=(0.8,10)\) case does not actually exceed the firing rate threshold \(\non=10\), even though the input pulse triggers MFEs. The \((b,J_{\max})=(0.8,20)\) case, on the other hand, triggers more intensive MFEs than the former case, which push the firing rate to exceed the firing rate threshold \(\non=10\). Nevertheless, it is shown that the multi-scale solver produces an accurate approximation to the microscopic results for both cases, indicating the reliability of the multi-scale solver.

Also, we mention that the multi-scale solver remarkably reduces the time cost of simulating the network and we present the system times of the two schemes in Table 3. 

\begin{center}
  \begin{tabular}{c|cccc}
   \toprule
   \((b,J_{\max})=(0.8,20)\) & \(L=10^4\) & \(L=15^4\) & \(L=20^4\) & \(L=25^4\)  \\
   \midrule
   microscopic  & 4.60 & 19.65 & 57.10  & 128.16\\
   multi-scale & 2.49 & 6.22 & 13.30& 27.25\\
   \midrule
   \((b,J_{\max})=(0.8,10)\) & \(L=10^4\) & \(L=15^4\) & \(L=20^4\) & \(L=25^4\)  \\
   \midrule
   microscopic  & 4.57 & 18.51 & 55.57  & 129.38\\
   multi-scale & 1.96 & 3.10 & 4.26 & 5.51\\
   \bottomrule
  \end{tabular}
 \captionof{table}{The table presents the system times to simulate a \(L\) neurons network with input current \eqref{T6I0} with network parameter \((b,J_{\max})\). }
 \label{TabT6}
 \end{center}

\section{Conclusion and perspectives}

The paper proposes a multi-scale solver for the noisy leaky integrate-and-fire (NLIF) networks that inherits the low cost of the macroscopic solver and the reliability of the microscopic solver. The multi-scale solver uses the macroscopic solver when the firing rate of the simulated network is low, while it switches to the microscopic solver when the firing rate tends to blow up. The validity of the multi-scale solver is analyzed from two perspectives: firstly, sufficient conditions are provided that guarantee the mean-field approximation of the macroscopic model and rigorous numerical analysis on simulation errors when coupling the two solvers; secondly, the numerical performance of the multi-scale solver is validated through simulating several large neuron networks with extensive experiments.

It is definitely worth investigating the extension of such a multi-scale method to a large network of both excitatory and inhibitory neurons. In this case, not only the structure of the MFEs are more complex, but also the integrate-and-fire model admits more diverse solutions. We save such explorations for future studies. 

\section*{Acknowledgements}

 The work of Z. Zhou is partially supported by the National Key R\&D Program of China (Project No. 2021YFA1001200, 2020YFA0712000), and the National Natural Science Foundation of China (Grant No. 12031013, 12171013). We thank Jiwei Zhang for the helpful discussions.

\section*{\textbf{Appendix.} The diffusion approximation}

Here we discuss the deduction of the diffusion approximation of the external synaptic currents \(I^\ite\) through the central limit theorem (CLT). We remark that the diffusion is obtained under the assumption that the network size \(L\gg 1\) and each individual neuron spikes independently.

Consider a NLIF neuron network with one species of neuron. We denote the collective synaptic currents of an individual neuron \(j\) as a stochastic process \(U_j(t)\) (\(j=1,2,\cdots,L\) ). Recall that we have assumed an individual neuron emit firing events with probability \(N(t)\) per unit time and each spike contributes an increment of \(J\) to the synaptic currents \(U_j(t)\). In other words, \(\frac{1}{J}U_j(t)\) is a nonhomogeneous Possion point process with time-dependent rate function \(N(t)\). Therefore, for a fixed \(t\), the mean and variance of \(I(t)\) follows
\beq 
\BE(U_j(t))=J\int_0^t N(s)\dd s,\, \textrm{var}(U_j(t))=J^2\int_0^t N(s)\dd s.
\eeq

For a \(L\) neurons network, the collective synaptic current \(U(t)\) follows
\beq 
U(t)=\sum_{j=1}^L U_j(t)
\eeq
With the condition that each neuron fires independently, we assume the stochastic processes \(U_j(t)\) are independent. As the network size \(L\to\infty\), we naturally adopt the CLT to study \(U(t)\):
\beq \label{Apd:limit}
\frac{1}{\sqrt{L}}\kl U(t)-JL\int_0^t N(s)\dd s\kr\overset{\textrm{d.}}{\longrightarrow}\CN\kl 0, J^2\int_0^t N(s)\dd s\kr,
\eeq
where \(\CN(\mu,\sigma^2)\) denotes a normal distribution. Recall that every \(U_j(t)\) have independent increments, we may approximate the stochastic fluctuations on the left hand side of Equation \eqref{Apd:limit} with an Ito process
\beq 
\frac{1}{\sqrt{L}}\kl U(t)-JL\int_0^t N(s)\dd s\kr\approx J\int_0^t\sqrt{N(s)}\dd B_s,
\eeq
where we have the Ito isometric \(\BE\left[\int_0^t\sqrt{N(s)}\dd B_s\right]^2=\int_0^t N(s)\dd s\). Finally, we arrive at the diffusion approximation to the total synaptic current \(U(t)\) when \(L\gg 1\): 
\beq 
\dd U(t)\approx JLN(t)\dd t+J\sqrt{LN(t)}\dd B_t.
\eeq
Recall that the external current \(I^\ite(t)\) models the rate of currents accepted by an individual neuron per unit time within the entire network, then
\beq 
\dd U(t)=I^{\ite}(t) \dd t\approx JLN(t)\dd t+J\sqrt{LN(t)}\dd B_t.
\eeq

\end{document}